\documentclass[12pt]{amsart}
\usepackage{epic,eepic}
\usepackage{amscd,amssymb}
 \newlength{\baseunit}               
\newcommand{\bpf}{\noindent {\em Proof. }}
\newcommand{\epf}{\qed \vspace{+10pt}}

\newtheorem{tm}{Theorem}[section]
\newtheorem{pr}[tm]{Proposition}
\newtheorem{lm}[tm]{Lemma}

\newcommand{\proj}{\mathbb P}

\newcommand{\eff}{\mathbb F}
\newcommand{\fn}{{\mathbb F}_n}
\newcommand{\kfn}{K_{\fn}}
\newcommand{\aff}{\mathbb A}

\newcommand{\oh}{{\mathcal{O}}}
\newcommand{\cc}{{\mathcal{C}}}

\newcommand{\cj}{{\mathcal{J}}}

\newcommand{\cs}{{\mathcal{S}}}

\newcommand{\zed}{{\mathbb{Z}}}
\newcommand{\tth}{\text{th}}
\newcommand{\ti}{\operatorname{irr}}
\newcommand{\tpi}{\tilde{\pi}}
\newcommand{\tc}{\tilde{C}}
\newcommand{\tg}{\tilde{g}}
\newcommand{\al}{\alpha}
\newcommand{\be}{\beta}
\newcommand{\ga}{\gamma}
\newcommand{\Ga}{\Gamma}
\newcommand{\de}{\delta}
\newcommand{\De}{\Delta}
\newcommand{\ka}{\kappa}
\newcommand{\si}{\sigma}
\newcommand{\Up}{\Upsilon}
\newcommand{\vdg}{V^{D,g}}
\newcommand{\ab}{(\al,\be)} 
\newcommand{\abG}{(\al,\be,\Ga)} 
\newcommand{\la}{\lambda}
\newcommand{\lcm}{\operatorname{lcm}}
\newcommand{\coker}{\operatorname{coker}}
\newcommand{\Sym}{\operatorname{Sym}}

\newcommand{\Def}{\operatorname{Def}}
\newcommand{\idim}{\operatorname{idim}}

\newcommand{\mbar}{\overline{M}}

\newcommand{\nm}{n \text{ mod }2}
\newcommand{\Nt}{N_{\operatorname{tors}}}

\begin{document}
\title{Counting Curves of Any Genus on Rational Ruled Surfaces}
\author{Ravi Vakil}
\date{\today}
\begin{abstract}
In this paper we study the geometry of the {\em Severi varieties}
parametrizing curves on the rational ruled surface $\fn$.  We compute the
number of such curves through the appropriate number of fixed general
points on $\fn$ (Theorem \ref{recursion}), and the number of such curves
which are irreducible (Theorem \ref{irecursion}).  These numbers are known
as {\em Severi degrees}; they are the degrees of unions of components of
the Hilbert scheme.  As (i) $\fn$ can be deformed to $\eff_{n+2}$, (ii)
the Gromov-Witten invariants are deformation-invariant, and (iii) the
Gromov-Witten invariants of $\eff_0$ and $\eff_1$ are enumerative, Theorem
\ref{irecursion} computes the genus $g$ Gromov-Witten invariants of all
$\fn$.  (The genus 0 case is well-known.)  The arguments are given in
sufficient generality to also count plane curves in the style of
L. Caporaso and J. Harris and to lay the groundwork for computing higher
genus Gromov-Witten invariants of blow-ups of the plane at up to five
points (in a future paper).
\end{abstract}
\maketitle
\pagestyle{plain}

\tableofcontents

\section{Introduction}

In this paper we study the geometry of the {\em Severi varieties}
parametrizing curves on the rational ruled surface $\fn = \proj
(\oh_{\proj^1} \oplus \oh_{\proj^1}(n))$ ($n \ge 0$) in a given divisor
class.  We compute the number of such curves through the appropriate number
of fixed general points on $\fn$ (Theorem \ref{recursion}), and the number
of such curves which are irreducible (Theorem \ref{irecursion}).  These
numbers are known as {\em Severi degrees}; they are the degrees of unions
of components of the Hilbert scheme.  As (i) $\fn$ can be deformed to
$\eff_{n+2}$, (ii) the Gromov-Witten invariants are deformation-invariants,
and (iii) the Gromov-Witten invariants of $\eff_0$ and $\eff_1$ are
enumerative, Theorem \ref{irecursion} computes the genus $g$ Gromov-Witten
invariants of all $\fn$ (Section \ref{gwenumerative}).  (The genus 0 case
is now well-known; it follows from the associativity of quantum cohomology.
See [KM], [DI] and [K1] for discussion.)  The arguments are given in
sufficient generality to also count plane curves in the style of
L. Caporaso and J. Harris (cf. [CH3]) and to lay the groundwork for
computing higher genus Gromov-Witten invariants of blow-ups of the plane at
up to five points (cf. [V3]).

Such a classical enumerative question has recently been the object of study
by many people.  Ideas from mathematical physics (cf. the inspiring [KM]
and [DI]) have yielded formulas when $g=0$ (via associativity relations in
quantum cohomology).  In March 1994, S. Kleiman and R. Piene found an
elegant recursive formula for the genus 0 Gromov-Witten invariants for all
$\fn$, and found empirically that many of the invariants of $\fn$ were
enumerative and the same as those for $\eff_{n+2}$, supporting the
conjecture (since proved) that quantum cohomology exists and is a
deformation invariant (see [K1] p. 22 for more information, and an
introduction to applications of quantum cohomology to enumerative
geometry).  Z. Ran solved the analogous (enumerative) problem for curves of
arbitrary genus on $\proj^2$ by degenerations methods (cf. [R]), and
Caporaso and Harris gave a second solution by different degeneration
methods (cf. [CH3]).  These numbers for irreducible curves are also the
genus $g$ Gromov-Witten invariants of $\proj^2$ (Section
\ref{gwenumerative}).  D. Abramovich and A. Bertram have used excess
intersection and the moduli space of stable maps to calculate generalized
Severi degrees for rational curves in all classes on $\eff_2$, and for
rational curves in certain classes on $\fn$ ([AB]).  (The author has used a
similar idea for curves of arbitrary genus in certain classes on $\fn$.)
In [CH1] and [CH2], Caporaso and Harris found recursive formulas for these
numbers when $g=0$ on $\eff_0$, $\eff_1$, $\eff_2$, and $\eff_3$, and on
certain classes on general $\fn$.  D. Coventry has also recently derived a
recursive formula for the number of rational curves in {\em any} class on
$\fn$ ([Co]), subsuming many earlier results.  In another direction,
extending work of I. Vainsencher ([Va]), Kleiman and Piene have examined
systems with an arbitrary, but fixed, number $\de$ of nodes ([K2]).  The
postulated number of $\de$-nodal curves is given (conjecturally) by a
polynomial, and they determine the leading coefficients, which are
polynomials in $\de$.  L. G\"{o}ttsche has recently conjectured a
surprisingly simple generating function ([G]) for these polynomials which
reproduce the results of Vainsencher as well as Kleiman and Piene and
experimentally reproduce the numbers of [CH3], [V3], S.T. Yau and
E. Zaslow's count of rational curves on K3-surfaces ([YZ]), and others.
The numbers of curves are expressed in terms of four universal power
series, two of which G\"{o}ttsche gives explicitly as quasimodular forms.

The philosophy here is that of Caporaso and Harris in [CH3]: we degenerate
the point conditions to lie on $E$ one at a time.  Our perspective,
however, is different: we use the moduli space of stable maps rather
than the Hilbert scheme.

The author is grateful to J. Harris for originally suggesting this problem
and for his contagious enthusiasm, and to D. Abramovich, E. Getzler, and T.
Graber for
fruitful discussions on several parts of the argument.  The exposition was
strengthened immeasurably thanks to advice from S. Kleiman.  D. Watabe
provided useful comments on earlier drafts.  L. G\"{o}ttsche's maple
program implementing the algorithm provided examples, and it is a pleasure
to acknowledge him here.

This research was supported (at different times)
by a NSERC 1967 Fellowship and a Sloan Dissertation Fellowship.  The bulk
of this paper was written at the Mittag-Leffler Institute, and the author
is grateful for the warmth and hospitality of the Institute staff.

\subsection{Background}
\label{background}
We work over the complex numbers.  Most of the arguments will be in
some generality, so that they can be invoked in [V3] to count curves
of arbitrary genus in any divisor class on the blow-up of the plane at
up to five points.

The Picard group of $\fn$ is $\zed^2$, with generators corresponding
to the fiber of the projective bundle $F$ and a section $E$ of
self-intersection $-n$; $E$ is unique if $n>0$.  Let $S$ be the class
$E+nF$.  (This class is usually denoted $C$, but we use nonstandard
notation to prevent confusion with the source of a map $(C,\pi)$.)
The canonical bundle $\kfn$ is $-(S+E+2F)$.  

Throughout this paper, $X$ will be $\fn$.  Unless otherwise explicitly
stated, we will use only the following properties of $(X,E)$.

\begin{enumerate}
\item[{\bf P1.}]  $X$ is a smooth surface and $E \cong \proj^1$ is a divisor on $X$.
\item[{\bf P2.}]
The surface $X \setminus E$ is minimal, i.e. contains no (-1)-curves.
\item[{\bf P3.}]
The divisor class $K_X + E$ is negative on every curve on $X$.
\item[{\bf P4.}] 
If $D$ is an effective divisor such that $-(K_X + E) \cdot D = 1$,
then $D$ is smooth.
\end{enumerate}

Property P2 could be removed by modifying the arguments very slightly,
but there seems to be no benefit of doing so.  Properties P3 and P4
would follow if $-(K_X + E)$ were very ample, which is true in all
cases of interest here.

Notice that if $L$ is a line on $\proj^2$ then $(X,E) = (\proj^2,L)$
also satisfies properties P1--P4.  The resulting formulas for
Severi degrees of $\proj^2$ are then those of [CH3].  Theorem
\ref{recursion} would become Theorem 1.1 of [CH3], and Theorem
\ref{irecursion} would give a recursive formula for irreducible genus
$g$ curves (which are the genus $g$ Gromov-Witten invariants).  

If $C$ is a smooth conic on $\proj^2$, then $(X,E) = (\proj^2,C)$ also
satisfies properties P1--P4.  This will be the basis of the computation of
higher genus Gromov-Witten invariants of blow-ups of $\proj^2$ at up
to 5 points in [V3].

For any sequence $\al = (\al_1, \al_2, \dots)$ of nonnegative integers
with all but finitely many $\al_i$ zero, set 
$$
| \al | = \al_1 + \al_2 + \al_3 + \dots 
$$
$$
I  \al
 = \al_1 + 2\al_2 + 3\al_3 + \dots 
$$
$$
I^\al
 = 1^{\al_1}  2^{\al_2}  3^{\al_3}  \dots 
$$
and
$$
\al ! = \al_1 ! \al_2 ! \al_3! \dots .
$$ We denote by $\lcm(\al)$ the least common multiple of the set $\#
\{ i : \al_i \neq 0 \}$.  The zero sequence will be denoted 0.

We denote by $e_k$ the sequence $(0, \dots, 0, 1, 0 , \dots)$ that is
zero except for a 1 in the $k^{\tth}$ term (so that any sequence $\al
= (\al_1, \al_2, \dots)$ is expressible as $\al = \sum \al_k e_k$).
By the inequality $\al \geq \al'$ we mean $\al_k \geq \al_k'$ for all
$k$; for such a pair of sequences we set $$
\binom \al {\al'} = {\frac { \al!}  {\al' ! (\al - \al')!}} = 
\binom {\al_1} {\al'_1}
\binom {\al_2} {\al'_2}
\binom {\al_3} {\al'_3} \dots .
$$
This notation follows [CH3].

For any divisor class $D$ on $X$, genus $g$, sequences $\al$ and $\be$, and
collections of points $\Ga = \{ p_{i,j} \}_{1 \leq j \leq \al_i}$ (not
necessarily distinct) of $E$ we define the {\em generalized Severi variety}
$\vdg(\al,\be,\Ga)$ to be the closure (in $|D|$) of the locus of reduced
curves $C$ in $X$ in divisor class $D$ of geometric genus $g$, not
containing $E$, with (informally) $\al_k$ ``assigned'' points of contact of
order $k$ and $\be_k$ ``unassigned'' points of contact of order $k$ with
$E$.  Formally, we require that, if $\nu: C^\nu \rightarrow C$ is the
normalization of $C$, then there exist $|\al|$ points $q_{i,j} \in C^\nu$,
$j=1, \dots, \al_i$ and $|\be|$ points $r_{i,j} \in C^\nu$, $j=1, \dots,
\be_i$ such that 
$$
 \nu(q_{i,j}) = p_{i,j} \quad \text{and} \quad \nu^*(E) = \sum i \cdot
q_{i,j} + \sum i \cdot r_{i,j}.
$$
If $I \al + I\be \neq D \cdot E$, $\vdg(\al,\be,\Ga)$ is empty.

For convenience, let 
$$
\Up = \Up^{D,g}(\be) := - (K_X + E) \cdot D +
|\be| + g-1.
$$
We will see that $\vdg\abG$ is a projective variety of pure dimension $\Up$
(Prop. \ref{idim}).  Let $N^{D,g}\abG$ be the number of points of
$\vdg\abG$ whose corresponding curve passes through $\Up$ fixed general
points of $X$.  Then $N^{D,g}\abG$ is the degree of the generalized Severi
variety (in the projective space $|D|$).  When the points $\{ p_{i,j}\}$
are distinct, we will see that $N\abG$ is independent of $\Ga$ (Section
\ref{recursivesection}); for simplicity we will then write $N\ab$.  The
main result of this paper is the following.
\begin{tm}
\label{recursion}
If $\dim \vdg(\al,\be)>0$, then 
\begin{eqnarray*}
N^{D,g}\ab = \sum_{\be_k > 0} k N^{D,g}(\al + e_k, \be-e_k) 
\\
+ \sum I^{\be'-\be} {\binom \al {\al'}} \binom {\be'}{\be} 
N^{D-E,g'}(\al',\be')
\end{eqnarray*}
where the second sum is taken over all $\al'$, $\be'$, $g'$ satisfying
$\al' \leq \al$, $\be' \geq \be$, $g-g' = |\be'-\be| - 1$, $I \al' + I \be' =
(D-E) \cdot E$.
\end{tm}
When $X=\fn$, the condition $\Up^{D,g}(\be) = \dim \vdg\ab>0$ is equivalent to $(D,g,\be)
\neq (k F, 1-k,0)$.  With the ``seed data''
$$
N^{kF,1-k}(\al,0) = 
\begin{cases}
1 & \text{if $\al = k e_1$,} \\
0 & \text{otherwise,}
\end{cases}
$$
this formula inductively counts curves of any
genus in any divisor class of $\fn$.

In order to understand generalized Severi varieties, we will analyze
certain moduli spaces of maps.  Let $\mbar_g(X,D)'$ be the moduli space of
maps $\pi: C \rightarrow X$ where $C$ is complete, reduced, and nodal, $(C,
\pi)$ has finite automorphism group, and $\pi_* [C] = D$.  (The curve $C$
is not required to be irreducible.)  If $C'$ is any connected component of
$C$, the map $(C', \pi)$ is stable.  The moduli space $\mbar_g(X,D)$ is a
union of connected components of $\mbar_g(X,D)'$.

Let $D$, $g$, $\al$, $\be$, $\Ga$ be
as in the definition of $\vdg(\al,\be,\Ga)$ above.  Define the {\em
generalized Severi variety of maps} $\vdg_m(\al,\be,\Ga)$ to be the
closure in $\mbar_g(X,D)'$ of points representing maps $(C,\pi)$ where
each component of $C$ maps birationally to its image in $X$, no
component maps to $E$, and $C$ has (informally) $\al_k$ ``assigned''
points of contact of order $k$ and $\be_k$ ``unassigned'' points of
contact of order $k$ with $E$.  Formally, we require that there exist
$|\al|$ smooth points $q_{i,j} \in C$, $j=1, \dots, \al_i$ and $|\be|$
smooth points $r_{i,j} \in C$, $j=1, \dots, \be_i$ such that 
$$
\pi(q_{i,j}) = p_{i,j} \quad \text{and} \quad \pi^*(E) = \sum i \cdot q_{i,j} +
\sum i \cdot r_{i,j}.
$$
As before, where the dependence on the points $p_{i,j}$ is not relevant
--- for example, in the discussions of the dimensions or degrees of
generalized Severi varieties --- we will suppress the $\Ga$.

There is a natural rational map from each component of
$\vdg\abG$ to $\vdg_m\abG$, and the 
dimension of the image will be $\Up$.  We will prove:

\begin{pr}
\label{idim}
The components of $\vdg_m\abG$ have dimension at most $\Up$,
and the union of those with dimension exactly $\Up$ is
the closure of the image of $\vdg\abG$ in $\vdg_m\abG$.
\end{pr}

(This will be an immediate consequence of Theorem \ref{bigdim}.)

Assume now that the $\{ p_{i,j} \}$ are distinct.
Fix $\Up$ general points $s_1$, \dots, $s_\Up$ on $X$.  The image of
the maps in $\vdg_m\ab$ whose images pass through these points are
reduced.  ({\em Proof:} Without loss of generality, restrict to the
union $V$ of those components of $\vdg_m\ab$ with dimension $\Up$.  By
Proposition \ref{idim}, the subvariety of $V$ corresponding to maps
whose images are {\em not} reduced contains no components of
$V$ and hence has dimension less than $\Up$.  Thus no image of such a
map passes through $s_1$, \dots, $s_{\Up}$.)

Therefore, if $H$ is the divisor class on $\vdg_m\ab$ corresponding to
requiring the image curve to pass through a fixed point of $X$, then
$$
N^{D,g}\ab = H^\Up.  
$$ 

Define the {\em intersection dimension} of a family $V$ of maps to
$X$ (denoted $\idim V$) as the maximum number $n$ of general
points $s_1$, \dots $s_n$ on $X$ such that there is a map $\pi: C
\rightarrow X$ in $V$ with $\{ s_1, \dots, s_n \} \subset
\pi(C)$.  Clearly $\idim V \leq \dim V$.

Our strategy is as follows.  Fix a general point $q$ of $E$.  Let
$H_q$ be the Weil divisor on $\vdg_m\abG$ corresponding to maps with images
containing $q$.  We will find the components of $\vdg_m\abG$ with
intersection dimension $\Up-1$ and relate them to
$V^{D',g'}_m(\al',\be', \Ga')$ for appropriately chosen $D'$, $g'$, $\al'$,
$\be'$, $\Ga'$.  Then we compute the multiplicity with which each of these
components appears.  Finally, we derive a
recursive formula for $N^{D,g}\ab$ (Theorem \ref{recursion}).

Analogous definitions can be made of spaces $W^{D,g}\abG$ and
$W^{D,g}_m\abG$ parametrizing irreducible curves.  The arguments in
this case are identical, resulting in a recursive formula for
$N_{\ti}^{D,g}\ab$, the number of irreducible genus $g$ curves in class
$D$ intersecting $E$ as determined by $\al$ and $\be$, passing through
$\Up$ fixed general points of $X$:
\begin{tm}
\label{irecursion}
If $\dim W^{D,g}\ab>0$, then 
\begin{eqnarray*}
N_{\ti}^{D,g}\ab &=&  \sum_{\beta_k > 0} k N_{\ti}^{D,g}(\alpha + e_k,
\beta - e_k)
\\
& & + \sum \frac 1 \si \binom {\al} {\al^1, \dots, \al^l, \al- \sum \al^i} \binom {\Up^{D,g}(\be)-1} {\Up^{D^1,g^1}(\be^1), \dots, \Up^{D^l,g^l}(\be^l)} \\
& & \cdot  \prod_{i=1}^l \binom {\be^i} {\ga^i} I^{\be^i - \ga^i} N_{\ti}^{D^i,g^i}(\al^i,\be^i) 
\end{eqnarray*}
where the second sum runs over choices of $D^i, g^i, \alpha^i, \beta^i,
\gamma^i$ ($1 \le i \le l$), where $D^i$ is a divisor class, $g^i$ is a
non-negative integer, $\alpha^i$, $\be^i$, $\ga^i$ are sequences of
non-negative integers, $\sum D^i = D-E$, $\sum \ga^i = \be$, $\be^i \gneq
\ga^i$, and $\si$ is the number of symmetries of the set $\{
(D^i,g^i,\al^i,\be^i,\ga^i) \}_{1 \leq i \leq l}$.
\end{tm}
In the second sum, for the summand to be non-zero, one must also have
$\sum \al^i \leq \al$, and $I \al^i + I \be^i =
D^i \cdot E$.  When $X=\fn$, the condition $\dim W^{D,g}\ab>0$ is
equivalent to $(D,\be) \neq (F, 0)$.  Thus with the ``seed data''
$N_{\ti}^{F,0}(e_1,0) = 1$,
 this formula inductively
counts irreducible curves of any genus in any divisor class of $\fn$.

\subsection{Examples}
As an example of the algorithm in action, we calculate $N^{4S,1}(0,0) =
225$ on $\eff_1$.  (This is also the number of two-nodal elliptic plane
quartics through 11 fixed general points.)  There are a finite number of
such elliptic curves through 11 fixed general points on $\eff_1$.  We
calculate the number by specializing the fixed points to lie on $E$ one at
a time, and following what happens to the finite number of curves.

After the first specialization, the curve must contain $E$ (as $4S \cdot E
= 0$, any representative of $4S$ containing a point of $E$ must contain
all of $E$).  The residual curve is in class $3S+F$.  Theorem
\ref{recursion} gives
$$
N^{4S,1}(0,0) = N^{3S+F,1}(0,e_1).
$$

After specializing a second point $q$ to lie on $E$, two things could
happen to the elliptic curve.  First, the limit curve could remain smooth,
and pass through the fixed point $q$ of $E$.  This will happen
$N^{3S+F,1}(e_1,0)$ times.  Second, the curve could contain $E$.  Then the
residual curve $C'$ has class $2S+2F$, and is a nodal curve intersecting
$E$ at two distinct points.  Of the two nodes of the original curve $C$,
one goes to the node of $C'$, and the other tends to one of the
intersection of $C'$ with $E$.  The choice of the two possible limits of
the node gives a multiplicity of 2.  Theorem \ref{recursion} gives
$$
N^{3S+F,1}(0,e_1) = N^{3S+F,1}(e_1,0) + 2 N^{2S+2F,1}(0,2e_1).
$$

Now $N^{2S+2F,1}(0,2e_1)$, the number of nodal curves in the linear system
$| 2S+2F|$, can be calculated to be 20 by further degenerations or by the
well-known calculation of the degree of the hypersurface of singular
sections in any linear system.  This calculation is omitted.

The number $N^{3S+F,1}(e_1,0)$ is calculated by specializing another point
to be a general point of $E$.  The limit curve will be of one of three
forms; in each case the limit must contain $E$, and the residual curve $C''$
is in the class $2S+2F$.
\begin{enumerate}
\item  The curve $C''$ could have geometric genus 0 and intersect $E$ at two
points.  There are two subcases:  $C''$ could be irreducible, or it could
consist of a fiber $F$ and a smooth elliptic curve in the class $2S+F$.
These cases happen $N^{2S+2F,0}(0,2e_1)$ times.
\item The curve $C''$ has geometric genus 1 and is tangent to $E$ at a
general point.  This happens $N^{2S+2F,1}(0,e_2)$ times.  Each of these
curves is the limit of {\em two} curves, so there is a multiplicity of 2.
(This multiplicity is not obvious.)
\item The curve $C''$ is smooth, and passes through the point $q \in E$.
This happens $N^{2S+2F,1}(e_1,e_1)$ times.  
\end{enumerate}

Theorem \ref{recursion} gives us
$$
N^{3S+F,1}(e_1,0) = N^{2S+2F,0}(0,2e_1) + 2 N^{2S+2F,1}(0,e_2) +
N^{2S+2F,1}(e_1,e_1).
$$

One can continue and calculate 
$$N^{2S+2F,0}(0,2e_1) = 105, \: N^{2S+2F,1}(0,e_2) = 30, \:
N^{2S+2F,1}(e_1,e_1)=20.
$$
Then we can
recursively calculate $N^{4S,1}(0,0)$:
\begin{eqnarray*}
N^{3S+F,1}(e_1,0) &=& N^{2S+2F,0}(0,2e_1) + 2 N^{2S+2F,1}(0,e_2) +
N^{2S+2F,1}(e_1,e_1) \\
&=& 105 + 2 \cdot 30 + 20 \\
&=& 185
\end{eqnarray*}
\begin{eqnarray*}
\text{so } N^{4S,1}(0,0) &=& N^{3S+F,1}(0,e_1) \\
 &=& N^{3S+F,1}(e_1,0) + 2 N^{2S+2F,1}(0,2e_1) \\
 &=& 185 + 2 \cdot 20 \\
 &=& 225.
\end{eqnarray*}

\begin{figure}
\begin{center}
	   \setlength{\unitlength}{.1\baseunit}
{%
\begingroup\makeatletter\ifx\SetFigFont\undefined
\def\x#1#2#3#4#5#6#7\relax{\def\x{#1#2#3#4#5#6}}%
\expandafter\x\fmtname xxxxxx\relax \def\y{splain}%
\ifx\x\y   
\gdef\SetFigFont#1#2#3{%
  \ifnum #1<17\tiny\else \ifnum #1<20\small\else
  \ifnum #1<24\normalsize\else \ifnum #1<29\large\else
  \ifnum #1<34\Large\else \ifnum #1<41\LARGE\else
     \huge\fi\fi\fi\fi\fi\fi
  \csname #3\endcsname}%
\else
\gdef\SetFigFont#1#2#3{\begingroup
  \count@#1\relax \ifnum 25<\count@\count@25\fi
  \def\x{\endgroup\@setsize\SetFigFont{#2pt}}%
  \expandafter\x
    \csname \romannumeral\the\count@ pt\expandafter\endcsname
    \csname @\romannumeral\the\count@ pt\endcsname
  \csname #3\endcsname}%
\fi
\fi\endgroup
\begin{picture}(11037,15348)(0,-10)
\thicklines
\dottedline{135}(1800,1821)(3900,1821)
\dottedline{135}(4800,1821)(6900,1821)
\path(5250,2121)(5250,771)
\path(2100,546)	(2112.093,591.278)
	(2124.028,635.601)
	(2135.810,678.974)
	(2147.441,721.407)
	(2158.926,762.906)
	(2170.267,803.479)
	(2181.467,843.132)
	(2192.532,881.873)
	(2214.264,956.649)
	(2235.491,1027.864)
	(2256.240,1095.577)
	(2276.539,1159.847)
	(2296.415,1220.732)
	(2315.895,1278.290)
	(2335.007,1332.579)
	(2353.779,1383.658)
	(2372.238,1431.586)
	(2390.411,1476.420)
	(2408.326,1518.219)
	(2426.010,1557.041)
	(2460.796,1625.990)
	(2494.990,1683.732)
	(2528.810,1730.734)
	(2562.476,1767.465)
	(2630.227,1811.976)
	(2700.000,1821.000)

\path(2700,1821)	(2772.748,1759.423)
	(2798.844,1694.793)
	(2809.760,1655.590)
	(2819.422,1612.395)
	(2827.954,1565.655)
	(2835.482,1515.822)
	(2842.131,1463.343)
	(2848.025,1408.669)
	(2853.289,1352.248)
	(2858.049,1294.531)
	(2862.430,1235.965)
	(2866.556,1177.001)
	(2870.553,1118.088)
	(2874.546,1059.675)
	(2878.659,1002.211)
	(2883.018,946.146)
	(2887.748,891.929)
	(2892.974,840.009)
	(2898.820,790.835)
	(2905.412,744.857)
	(2912.875,702.525)
	(2921.334,664.287)
	(2941.739,601.891)
	(3000.000,546.000)

\path(3000,546)	(3068.779,563.740)
	(3133.580,625.672)
	(3164.970,674.720)
	(3195.942,736.629)
	(3211.332,772.596)
	(3226.690,812.004)
	(3242.039,854.929)
	(3257.404,901.448)
	(3272.809,951.634)
	(3288.277,1005.565)
	(3303.834,1063.315)
	(3319.502,1124.960)
	(3335.306,1190.575)
	(3351.271,1260.236)
	(3367.419,1334.020)
	(3375.570,1372.480)
	(3383.775,1412.000)
	(3392.039,1452.587)
	(3400.364,1494.252)
	(3408.753,1537.005)
	(3417.208,1580.853)
	(3425.734,1625.808)
	(3434.333,1671.878)
	(3443.008,1719.073)
	(3451.762,1767.402)
	(3460.598,1816.875)
	(3469.519,1867.500)
	(3478.528,1919.289)
	(3487.629,1972.249)
	(3496.823,2026.391)
	(3506.114,2081.724)
	(3515.506,2138.257)
	(3525.000,2196.000)

\path(4800,996)	(4852.704,1004.433)
	(4903.528,1012.790)
	(4952.515,1021.083)
	(4999.711,1029.328)
	(5045.159,1037.537)
	(5088.906,1045.725)
	(5130.994,1053.904)
	(5171.470,1062.090)
	(5210.377,1070.295)
	(5247.760,1078.534)
	(5318.134,1095.165)
	(5382.948,1112.095)
	(5442.559,1129.432)
	(5497.324,1147.288)
	(5547.601,1165.771)
	(5593.747,1184.991)
	(5636.118,1205.059)
	(5710.966,1248.176)
	(5775.000,1296.000)

\path(5775,1296)	(5818.483,1342.054)
	(5857.892,1401.268)
	(5893.796,1471.381)
	(5910.611,1509.817)
	(5926.763,1550.131)
	(5942.322,1592.040)
	(5957.360,1635.260)
	(5971.948,1679.509)
	(5986.156,1724.505)
	(6000.056,1769.964)
	(6013.719,1815.605)
	(6027.215,1861.145)
	(6040.616,1906.301)
	(6053.993,1950.791)
	(6067.416,1994.331)
	(6080.957,2036.640)
	(6094.687,2077.435)
	(6108.676,2116.433)
	(6122.996,2153.352)
	(6152.911,2219.820)
	(6185.002,2274.580)
	(6219.835,2315.371)
	(6257.978,2339.931)
	(6300.000,2346.000)

\path(6300,2346)	(6370.180,2324.595)
	(6431.504,2271.241)
	(6459.188,2231.416)
	(6485.070,2182.203)
	(6509.288,2123.134)
	(6531.979,2053.744)
	(6542.794,2015.032)
	(6553.279,1973.564)
	(6563.451,1929.283)
	(6573.327,1882.129)
	(6582.924,1832.045)
	(6592.260,1778.972)
	(6601.351,1722.851)
	(6610.215,1663.624)
	(6618.868,1601.233)
	(6627.329,1535.620)
	(6635.613,1466.726)
	(6643.739,1394.493)
	(6647.748,1357.106)
	(6651.724,1318.862)
	(6655.668,1279.754)
	(6659.584,1239.775)
	(6663.473,1198.917)
	(6667.337,1157.174)
	(6671.179,1114.537)
	(6675.000,1071.000)

\put(2865,9921){\blacken\ellipse{150}{150}}
\put(2865,9921){\ellipse{150}{150}}
\put(10081,7228){\blacken\ellipse{150}{150}}
\put(10081,7228){\ellipse{150}{150}}
\put(4283,4521){\blacken\ellipse{150}{150}}
\put(4283,4521){\ellipse{150}{150}}
\put(6480,4528){\blacken\ellipse{150}{150}}
\put(6480,4528){\ellipse{150}{150}}
\put(7530,4521){\blacken\ellipse{150}{150}}
\put(7530,4521){\ellipse{150}{150}}
\dottedline{135}(3900,15321)(6000,15321)
\dottedline{135}(3900,12621)(6000,12621)
\dottedline{135}(2250,9921)(4350,9921)
\dottedline{135}(5400,9921)(7500,9921)
\dottedline{135}(300,7221)(2400,7221)
\dottedline{135}(3000,7221)(5100,7221)
\path(3450,7821)(3450,6471)
\dottedline{135}(5700,7221)(7800,7221)
\dottedline{135}(8700,7221)(10800,7221)
\dottedline{135}(3300,4521)(5400,4521)
\dottedline{135}(6150,4521)(8250,4521)
\dottedline{135}(8925,4521)(11025,4521)
\path(4875,13746)(4875,12996)
\path(4875,13746)(4875,12996)
\path(4845.000,13116.000)(4875.000,12996.000)(4905.000,13116.000)
\path(4725,11046)(3300,10071)
\path(4725,11046)(3300,10071)
\path(3382.096,10163.521)(3300.000,10071.000)(3415.977,10114.003)
\path(4950,11046)(5850,10146)
\path(4950,11046)(5850,10146)
\path(5743.934,10209.640)(5850.000,10146.000)(5786.360,10252.066)
\path(7425,8346)(9075,7596)
\path(7425,8346)(9075,7596)
\path(8953.342,7618.345)(9075.000,7596.000)(8978.170,7672.967)
\path(2700,8346)(1200,7671)
\path(2700,8346)(1200,7671)
\path(1297.120,7747.601)(1200.000,7671.000)(1321.742,7692.886)
\path(3300,8346)(3375,7896)
\path(3300,8346)(3375,7896)
\path(3325.680,8009.435)(3375.000,7896.000)(3384.864,8019.299)
\path(3900,8346)(5925,7371)
\path(3900,8346)(5925,7371)
\path(5803.865,7396.028)(5925.000,7371.000)(5829.894,7450.088)
\path(4500,8346)(8475,7446)
\path(4500,8346)(8475,7446)
\path(8351.338,7443.240)(8475.000,7446.000)(8364.587,7501.758)
\path(6600,5646)(4650,4746)
\path(6600,5646)(4650,4746)
\path(4746.383,4823.526)(4650.000,4746.000)(4771.527,4769.048)
\path(9150,5571)(7575,4896)
\path(9150,5571)(7575,4896)
\path(7673.480,4970.845)(7575.000,4896.000)(7697.115,4915.696)
\path(9525,5571)(9525,5196)
\path(9525,5571)(9525,5196)
\path(9495.000,5316.000)(9525.000,5196.000)(9555.000,5316.000)
\path(4125,2946)(3600,2346)
\path(4125,2946)(3600,2346)
\path(3656.443,2456.064)(3600.000,2346.000)(3701.598,2416.554)
\path(4425,2946)(5175,2346)
\path(4425,2946)(5175,2346)
\path(5062.555,2397.537)(5175.000,2346.000)(5100.037,2444.389)
\path(4125,15096)	(4173.021,15088.102)
	(4219.346,15080.480)
	(4264.016,15073.127)
	(4307.074,15066.036)
	(4348.560,15059.200)
	(4388.514,15052.612)
	(4426.979,15046.266)
	(4463.996,15040.154)
	(4533.849,15028.607)
	(4598.404,15017.916)
	(4657.988,15008.026)
	(4712.932,14998.882)
	(4763.567,14990.430)
	(4810.220,14982.613)
	(4853.222,14975.378)
	(4892.903,14968.670)
	(4963.617,14956.612)
	(5025.000,14946.000)

\path(5025,14946)	(5063.401,14936.907)
	(5107.472,14922.657)
	(5156.372,14904.345)
	(5209.263,14883.071)
	(5265.307,14859.930)
	(5323.666,14836.020)
	(5383.499,14812.438)
	(5443.969,14790.281)
	(5504.236,14770.647)
	(5563.463,14754.633)
	(5620.811,14743.335)
	(5675.440,14737.851)
	(5726.512,14739.279)
	(5773.189,14748.714)
	(5814.631,14767.256)
	(5850.000,14796.000)

\path(5850,14796)	(5879.378,14862.562)
	(5881.163,14904.179)
	(5874.810,14948.033)
	(5860.669,14991.632)
	(5839.087,15032.486)
	(5775.000,15096.000)

\path(5775,15096)	(5716.131,15124.781)
	(5655.798,15142.307)
	(5594.390,15149.767)
	(5532.294,15148.348)
	(5469.900,15139.237)
	(5407.595,15123.621)
	(5345.770,15102.690)
	(5284.811,15077.629)
	(5225.109,15049.626)
	(5167.051,15019.870)
	(5111.026,14989.547)
	(5057.423,14959.845)
	(5006.630,14931.952)
	(4959.036,14907.055)
	(4915.030,14886.342)
	(4875.000,14871.000)

\path(4875,14871)	(4800.684,14836.037)
	(4758.083,14810.208)
	(4712.541,14780.471)
	(4664.577,14748.068)
	(4614.708,14714.241)
	(4563.450,14680.233)
	(4511.321,14647.286)
	(4458.839,14616.643)
	(4406.520,14589.546)
	(4354.881,14567.237)
	(4304.441,14550.959)
	(4255.715,14541.953)
	(4209.222,14541.464)
	(4165.478,14550.732)
	(4125.000,14571.000)

\path(4125,14571)	(4081.429,14630.966)
	(4070.234,14672.235)
	(4066.260,14717.014)
	(4069.629,14762.291)
	(4080.461,14805.055)
	(4125.000,14871.000)

\path(4125,14871)	(4173.647,14900.907)
	(4226.344,14919.008)
	(4282.438,14926.576)
	(4341.281,14924.886)
	(4402.222,14915.211)
	(4464.609,14898.824)
	(4527.794,14876.998)
	(4591.125,14851.009)
	(4653.952,14822.128)
	(4715.625,14791.630)
	(4775.493,14760.789)
	(4832.906,14730.877)
	(4887.214,14703.169)
	(4937.766,14678.937)
	(4983.911,14659.457)
	(5025.000,14646.000)

\path(5025,14646)	(5065.411,14636.250)
	(5114.876,14625.528)
	(5171.773,14613.919)
	(5234.480,14601.509)
	(5301.377,14588.382)
	(5370.842,14574.624)
	(5441.253,14560.319)
	(5510.989,14545.552)
	(5578.428,14530.409)
	(5641.950,14514.974)
	(5699.932,14499.333)
	(5750.753,14483.570)
	(5824.427,14452.018)
	(5850.000,14421.000)

\path(5850,14421)	(5817.470,14393.827)
	(5781.578,14383.685)
	(5734.870,14375.524)
	(5679.048,14369.108)
	(5615.812,14364.199)
	(5546.864,14360.561)
	(5473.905,14357.959)
	(5436.453,14356.971)
	(5398.637,14356.154)
	(5360.668,14355.477)
	(5322.760,14354.911)
	(5247.976,14353.993)
	(5175.987,14353.164)
	(5108.494,14352.186)
	(5047.197,14350.824)
	(4993.799,14348.841)
	(4950.000,14346.000)

\path(4950,14346)	(4893.660,14341.052)
	(4828.787,14335.278)
	(4792.396,14332.013)
	(4752.964,14328.460)
	(4710.191,14324.590)
	(4663.774,14320.376)
	(4613.411,14315.792)
	(4558.799,14310.808)
	(4499.637,14305.399)
	(4435.623,14299.536)
	(4366.455,14293.192)
	(4291.829,14286.340)
	(4252.376,14282.714)
	(4211.445,14278.952)
	(4168.999,14275.048)
	(4125.000,14271.000)

\path(5475,9471)	(5513.531,9520.056)
	(5551.095,9566.998)
	(5587.738,9611.859)
	(5623.504,9654.668)
	(5658.436,9695.458)
	(5692.581,9734.257)
	(5758.686,9806.012)
	(5822.175,9870.181)
	(5883.404,9927.009)
	(5942.731,9976.745)
	(6000.514,10019.636)
	(6057.108,10055.929)
	(6112.872,10085.871)
	(6168.162,10109.710)
	(6223.335,10127.692)
	(6278.748,10140.065)
	(6334.758,10147.076)
	(6391.723,10148.972)
	(6450.000,10146.000)

\path(6450,10146)	(6494.995,10139.015)
	(6538.468,10126.576)
	(6580.479,10109.077)
	(6621.084,10086.912)
	(6660.343,10060.474)
	(6698.314,10030.159)
	(6770.627,9959.467)
	(6805.085,9919.880)
	(6838.489,9877.990)
	(6870.898,9834.190)
	(6902.370,9788.876)
	(6932.962,9742.440)
	(6962.735,9695.277)
	(6991.745,9647.781)
	(7020.053,9600.345)
	(7047.715,9553.363)
	(7074.791,9507.230)
	(7101.338,9462.339)
	(7127.416,9419.084)
	(7153.083,9377.858)
	(7178.398,9339.057)
	(7228.202,9270.300)
	(7277.297,9215.965)
	(7326.150,9179.203)
	(7375.228,9163.164)
	(7425.000,9171.000)

\path(7425,9171)	(7454.407,9198.393)
	(7476.600,9248.546)
	(7490.984,9314.202)
	(7496.966,9388.103)
	(7493.953,9462.990)
	(7481.350,9531.608)
	(7458.563,9586.697)
	(7425.000,9621.000)

\path(7425,9621)	(7381.010,9636.301)
	(7334.352,9636.711)
	(7285.381,9624.048)
	(7234.454,9600.127)
	(7181.924,9566.763)
	(7128.149,9525.774)
	(7073.482,9478.975)
	(7018.279,9428.182)
	(6962.895,9375.212)
	(6907.687,9321.881)
	(6853.008,9270.004)
	(6799.215,9221.397)
	(6746.662,9177.878)
	(6695.705,9141.261)
	(6646.699,9113.363)
	(6600.000,9096.000)

\path(6600,9096)	(6562.475,9087.528)
	(6522.998,9080.683)
	(6481.184,9075.492)
	(6436.649,9071.982)
	(6389.008,9070.180)
	(6337.877,9070.115)
	(6282.871,9071.814)
	(6223.605,9075.304)
	(6159.696,9080.612)
	(6090.758,9087.767)
	(6016.407,9096.795)
	(5977.081,9102.020)
	(5936.258,9107.724)
	(5893.890,9113.910)
	(5849.928,9120.582)
	(5804.325,9127.742)
	(5757.031,9135.395)
	(5708.000,9143.544)
	(5657.183,9152.192)
	(5604.533,9161.343)
	(5550.000,9171.000)

\path(2175,6921)	(2135.635,6866.620)
	(2098.041,6817.007)
	(2062.026,6772.023)
	(2027.397,6731.531)
	(1961.529,6663.473)
	(1898.899,6611.734)
	(1837.969,6575.216)
	(1777.201,6552.820)
	(1715.058,6543.447)
	(1650.000,6546.000)

\path(1650,6546)	(1603.592,6554.652)
	(1560.036,6568.709)
	(1519.167,6587.786)
	(1480.815,6611.499)
	(1410.996,6671.294)
	(1349.243,6745.018)
	(1320.972,6786.141)
	(1294.215,6829.592)
	(1268.805,6874.986)
	(1244.575,6921.939)
	(1221.357,6970.066)
	(1198.983,7018.983)
	(1177.288,7068.304)
	(1156.102,7117.646)
	(1135.260,7166.624)
	(1114.593,7214.852)
	(1093.935,7261.946)
	(1073.118,7307.523)
	(1051.974,7351.196)
	(1030.336,7392.582)
	(984.911,7466.952)
	(935.504,7527.556)
	(880.775,7571.316)
	(819.387,7595.157)
	(750.000,7596.000)

\path(750,7596)	(698.415,7582.407)
	(648.115,7558.728)
	(599.629,7525.907)
	(553.487,7484.885)
	(510.216,7436.605)
	(470.346,7382.010)
	(434.405,7322.041)
	(402.923,7257.641)
	(376.427,7189.753)
	(355.446,7119.317)
	(340.510,7047.278)
	(332.147,6974.577)
	(330.886,6902.156)
	(337.255,6830.958)
	(351.783,6761.925)
	(375.000,6696.000)

\path(375,6696)	(421.514,6621.488)
	(485.143,6564.824)
	(521.972,6542.467)
	(561.398,6523.711)
	(602.858,6508.268)
	(645.791,6495.851)
	(689.638,6486.174)
	(733.836,6478.948)
	(777.825,6473.887)
	(821.044,6470.703)
	(862.932,6469.110)
	(902.928,6468.820)
	(975.000,6471.000)

\path(975,6471)	(1018.304,6477.499)
	(1062.261,6491.707)
	(1106.831,6512.455)
	(1151.978,6538.575)
	(1197.665,6568.901)
	(1243.852,6602.264)
	(1290.503,6637.496)
	(1337.580,6673.429)
	(1385.045,6708.895)
	(1432.860,6742.728)
	(1480.988,6773.758)
	(1529.392,6800.818)
	(1578.032,6822.740)
	(1626.872,6838.356)
	(1675.874,6846.499)
	(1725.000,6846.000)

\path(1725,6846)	(1766.541,6838.423)
	(1807.409,6824.055)
	(1848.702,6802.016)
	(1891.519,6771.427)
	(1936.958,6731.411)
	(1986.120,6681.086)
	(2040.100,6619.576)
	(2100.000,6546.000)

\path(3075,7596)	(3140.607,7585.208)
	(3204.839,7574.532)
	(3267.708,7563.968)
	(3329.224,7553.514)
	(3389.397,7543.165)
	(3448.239,7532.919)
	(3505.761,7522.771)
	(3561.972,7512.718)
	(3616.884,7502.758)
	(3670.507,7492.885)
	(3722.852,7483.098)
	(3773.931,7473.393)
	(3823.752,7463.765)
	(3872.329,7454.213)
	(3919.670,7444.732)
	(3965.787,7435.318)
	(4010.690,7425.969)
	(4054.391,7416.681)
	(4096.900,7407.451)
	(4138.227,7398.275)
	(4178.384,7389.149)
	(4217.381,7380.071)
	(4291.938,7362.042)
	(4361.985,7344.162)
	(4427.606,7326.401)
	(4488.889,7308.734)
	(4545.919,7291.132)
	(4598.781,7273.569)
	(4647.562,7256.016)
	(4692.347,7238.446)
	(4733.222,7220.831)
	(4803.586,7185.359)
	(4859.341,7149.381)
	(4929.770,7075.023)
	(4945.816,7036.204)
	(4950.000,6996.000)

\path(4950,6996)	(4942.997,6957.032)
	(4924.507,6920.278)
	(4850.319,6852.757)
	(4793.248,6821.659)
	(4721.944,6792.117)
	(4680.740,6777.877)
	(4635.720,6763.965)
	(4586.798,6750.360)
	(4533.889,6737.040)
	(4476.906,6723.986)
	(4415.764,6711.176)
	(4350.378,6698.590)
	(4280.660,6686.208)
	(4206.526,6674.009)
	(4167.775,6667.972)
	(4127.889,6661.972)
	(4086.855,6656.009)
	(4044.664,6650.077)
	(4001.304,6644.177)
	(3956.764,6638.303)
	(3911.035,6632.455)
	(3864.105,6626.630)
	(3815.963,6620.825)
	(3766.600,6615.036)
	(3716.003,6609.263)
	(3664.163,6603.502)
	(3611.068,6597.751)
	(3556.709,6592.007)
	(3501.073,6586.267)
	(3444.151,6580.529)
	(3385.932,6574.791)
	(3326.404,6569.050)
	(3265.558,6563.302)
	(3203.383,6557.547)
	(3139.867,6551.780)
	(3075.000,6546.000)

\path(5700,6546)	(5738.531,6595.056)
	(5776.095,6641.998)
	(5812.738,6686.859)
	(5848.504,6729.668)
	(5883.436,6770.458)
	(5917.581,6809.257)
	(5983.686,6881.012)
	(6047.175,6945.181)
	(6108.404,7002.009)
	(6167.731,7051.745)
	(6225.514,7094.636)
	(6282.108,7130.929)
	(6337.872,7160.871)
	(6393.162,7184.710)
	(6448.335,7202.692)
	(6503.748,7215.065)
	(6559.758,7222.076)
	(6616.723,7223.972)
	(6675.000,7221.000)

\path(6675,7221)	(6719.995,7214.015)
	(6763.468,7201.576)
	(6805.479,7184.077)
	(6846.084,7161.912)
	(6885.343,7135.474)
	(6923.314,7105.159)
	(6995.627,7034.468)
	(7030.085,6994.880)
	(7063.489,6952.990)
	(7095.898,6909.190)
	(7127.370,6863.876)
	(7157.962,6817.440)
	(7187.735,6770.277)
	(7216.745,6722.781)
	(7245.053,6675.345)
	(7272.715,6628.363)
	(7299.791,6582.230)
	(7326.338,6537.339)
	(7352.416,6494.084)
	(7378.083,6452.858)
	(7403.398,6414.057)
	(7453.202,6345.300)
	(7502.297,6290.965)
	(7551.150,6254.203)
	(7600.228,6238.164)
	(7650.000,6246.000)

\path(7650,6246)	(7679.407,6273.393)
	(7701.600,6323.546)
	(7715.984,6389.202)
	(7721.966,6463.103)
	(7718.953,6537.990)
	(7706.350,6606.608)
	(7683.563,6661.697)
	(7650.000,6696.000)

\path(7650,6696)	(7606.010,6711.301)
	(7559.352,6711.711)
	(7510.381,6699.048)
	(7459.454,6675.127)
	(7406.924,6641.763)
	(7353.149,6600.774)
	(7298.482,6553.975)
	(7243.279,6503.183)
	(7187.895,6450.212)
	(7132.687,6396.881)
	(7078.008,6345.004)
	(7024.215,6296.397)
	(6971.662,6252.878)
	(6920.705,6216.261)
	(6871.699,6188.363)
	(6825.000,6171.000)

\path(6825,6171)	(6787.475,6162.528)
	(6747.998,6155.683)
	(6706.184,6150.492)
	(6661.649,6146.982)
	(6614.008,6145.180)
	(6562.877,6145.115)
	(6507.871,6146.814)
	(6448.605,6150.304)
	(6384.696,6155.612)
	(6315.758,6162.767)
	(6241.407,6171.795)
	(6202.081,6177.020)
	(6161.258,6182.724)
	(6118.890,6188.910)
	(6074.928,6195.582)
	(6029.325,6202.742)
	(5982.031,6210.395)
	(5933.000,6218.544)
	(5882.183,6227.192)
	(5829.533,6236.343)
	(5775.000,6246.000)

\path(8700,6846)	(8738.531,6895.056)
	(8776.095,6941.998)
	(8812.738,6986.859)
	(8848.504,7029.668)
	(8883.436,7070.458)
	(8917.581,7109.257)
	(8983.686,7181.012)
	(9047.175,7245.181)
	(9108.404,7302.009)
	(9167.731,7351.745)
	(9225.514,7394.636)
	(9282.108,7430.929)
	(9337.872,7460.871)
	(9393.162,7484.710)
	(9448.335,7502.692)
	(9503.748,7515.065)
	(9559.758,7522.076)
	(9616.723,7523.972)
	(9675.000,7521.000)

\path(9675,7521)	(9719.995,7514.015)
	(9763.468,7501.576)
	(9805.479,7484.077)
	(9846.084,7461.912)
	(9885.343,7435.474)
	(9923.314,7405.159)
	(9995.627,7334.468)
	(10030.085,7294.880)
	(10063.489,7252.990)
	(10095.898,7209.190)
	(10127.370,7163.876)
	(10157.962,7117.440)
	(10187.735,7070.277)
	(10216.745,7022.781)
	(10245.053,6975.345)
	(10272.715,6928.363)
	(10299.791,6882.230)
	(10326.338,6837.339)
	(10352.416,6794.084)
	(10378.083,6752.858)
	(10403.398,6714.057)
	(10453.202,6645.300)
	(10502.297,6590.965)
	(10551.150,6554.203)
	(10600.228,6538.164)
	(10650.000,6546.000)

\path(10650,6546)	(10679.407,6573.393)
	(10701.600,6623.546)
	(10715.984,6689.202)
	(10721.966,6763.103)
	(10718.953,6837.990)
	(10706.350,6906.608)
	(10683.563,6961.697)
	(10650.000,6996.000)

\path(10650,6996)	(10606.010,7011.301)
	(10559.352,7011.711)
	(10510.381,6999.048)
	(10459.454,6975.127)
	(10406.924,6941.763)
	(10353.149,6900.774)
	(10298.482,6853.975)
	(10243.279,6803.183)
	(10187.895,6750.212)
	(10132.687,6696.881)
	(10078.008,6645.004)
	(10024.215,6596.397)
	(9971.662,6552.878)
	(9920.705,6516.261)
	(9871.699,6488.363)
	(9825.000,6471.000)

\path(9825,6471)	(9787.475,6462.528)
	(9747.998,6455.683)
	(9706.184,6450.492)
	(9661.649,6446.982)
	(9614.008,6445.180)
	(9562.877,6445.115)
	(9507.871,6446.814)
	(9448.605,6450.304)
	(9384.696,6455.612)
	(9315.758,6462.767)
	(9241.407,6471.795)
	(9202.081,6477.020)
	(9161.258,6482.724)
	(9118.890,6488.910)
	(9074.928,6495.582)
	(9029.325,6502.742)
	(8982.031,6510.395)
	(8933.000,6518.544)
	(8882.183,6527.192)
	(8829.533,6536.343)
	(8775.000,6546.000)

\path(3375,3846)	(3413.531,3895.056)
	(3451.095,3941.998)
	(3487.738,3986.859)
	(3523.504,4029.668)
	(3558.436,4070.458)
	(3592.581,4109.257)
	(3658.686,4181.012)
	(3722.175,4245.181)
	(3783.404,4302.009)
	(3842.731,4351.745)
	(3900.514,4394.636)
	(3957.108,4430.929)
	(4012.872,4460.871)
	(4068.162,4484.710)
	(4123.335,4502.692)
	(4178.748,4515.065)
	(4234.758,4522.076)
	(4291.723,4523.972)
	(4350.000,4521.000)

\path(4350,4521)	(4394.995,4514.015)
	(4438.468,4501.576)
	(4480.479,4484.077)
	(4521.084,4461.912)
	(4560.343,4435.474)
	(4598.314,4405.159)
	(4670.627,4334.468)
	(4705.085,4294.880)
	(4738.489,4252.990)
	(4770.898,4209.190)
	(4802.370,4163.876)
	(4832.962,4117.440)
	(4862.735,4070.277)
	(4891.745,4022.781)
	(4920.052,3975.345)
	(4947.715,3928.363)
	(4974.791,3882.230)
	(5001.338,3837.339)
	(5027.416,3794.084)
	(5053.083,3752.858)
	(5078.398,3714.057)
	(5128.202,3645.300)
	(5177.297,3590.965)
	(5226.150,3554.203)
	(5275.228,3538.164)
	(5325.000,3546.000)

\path(5325,3546)	(5354.407,3573.393)
	(5376.600,3623.546)
	(5390.984,3689.202)
	(5396.966,3763.102)
	(5393.953,3837.990)
	(5381.350,3906.608)
	(5358.563,3961.697)
	(5325.000,3996.000)

\path(5325,3996)	(5281.010,4011.301)
	(5234.352,4011.711)
	(5185.381,3999.048)
	(5134.454,3975.127)
	(5081.924,3941.763)
	(5028.149,3900.774)
	(4973.482,3853.975)
	(4918.279,3803.183)
	(4862.895,3750.212)
	(4807.687,3696.881)
	(4753.008,3645.004)
	(4699.215,3596.397)
	(4646.662,3552.878)
	(4595.705,3516.261)
	(4546.699,3488.363)
	(4500.000,3471.000)

\path(4500,3471)	(4462.475,3462.528)
	(4422.998,3455.683)
	(4381.184,3450.492)
	(4336.649,3446.982)
	(4289.008,3445.180)
	(4237.877,3445.115)
	(4182.871,3446.814)
	(4123.605,3450.304)
	(4059.696,3455.612)
	(3990.758,3462.767)
	(3916.407,3471.795)
	(3877.081,3477.020)
	(3836.258,3482.724)
	(3793.890,3488.910)
	(3749.928,3495.582)
	(3704.325,3502.742)
	(3657.031,3510.395)
	(3608.000,3518.544)
	(3557.183,3527.192)
	(3504.533,3536.343)
	(3450.000,3546.000)

\path(6150,4146)	(6188.531,4195.056)
	(6226.095,4241.998)
	(6262.738,4286.859)
	(6298.504,4329.668)
	(6333.436,4370.458)
	(6367.581,4409.257)
	(6433.686,4481.012)
	(6497.175,4545.181)
	(6558.404,4602.009)
	(6617.731,4651.745)
	(6675.514,4694.636)
	(6732.108,4730.929)
	(6787.872,4760.871)
	(6843.162,4784.710)
	(6898.335,4802.692)
	(6953.748,4815.065)
	(7009.758,4822.076)
	(7066.723,4823.972)
	(7125.000,4821.000)

\path(7125,4821)	(7169.995,4814.015)
	(7213.468,4801.576)
	(7255.479,4784.077)
	(7296.084,4761.912)
	(7335.343,4735.474)
	(7373.314,4705.159)
	(7445.627,4634.468)
	(7480.085,4594.880)
	(7513.489,4552.990)
	(7545.898,4509.190)
	(7577.370,4463.876)
	(7607.962,4417.440)
	(7637.735,4370.277)
	(7666.745,4322.781)
	(7695.053,4275.345)
	(7722.715,4228.363)
	(7749.791,4182.230)
	(7776.338,4137.339)
	(7802.416,4094.084)
	(7828.083,4052.858)
	(7853.398,4014.057)
	(7903.202,3945.300)
	(7952.297,3890.965)
	(8001.150,3854.203)
	(8050.228,3838.164)
	(8100.000,3846.000)

\path(8100,3846)	(8129.407,3873.393)
	(8151.600,3923.546)
	(8165.984,3989.202)
	(8171.966,4063.102)
	(8168.953,4137.990)
	(8156.350,4206.608)
	(8133.563,4261.697)
	(8100.000,4296.000)

\path(8100,4296)	(8056.010,4311.301)
	(8009.352,4311.711)
	(7960.381,4299.048)
	(7909.454,4275.127)
	(7856.924,4241.763)
	(7803.149,4200.774)
	(7748.482,4153.975)
	(7693.279,4103.183)
	(7637.895,4050.212)
	(7582.687,3996.881)
	(7528.008,3945.004)
	(7474.215,3896.397)
	(7421.662,3852.878)
	(7370.705,3816.261)
	(7321.699,3788.363)
	(7275.000,3771.000)

\path(7275,3771)	(7237.475,3762.528)
	(7197.998,3755.683)
	(7156.184,3750.492)
	(7111.649,3746.982)
	(7064.008,3745.180)
	(7012.877,3745.115)
	(6957.871,3746.814)
	(6898.605,3750.304)
	(6834.696,3755.612)
	(6765.758,3762.767)
	(6691.407,3771.795)
	(6652.081,3777.020)
	(6611.258,3782.724)
	(6568.890,3788.910)
	(6524.928,3795.582)
	(6479.325,3802.742)
	(6432.031,3810.395)
	(6383.000,3818.544)
	(6332.183,3827.192)
	(6279.533,3836.343)
	(6225.000,3846.000)

\path(9075,3771)	(9087.093,3816.278)
	(9099.028,3860.601)
	(9110.810,3903.974)
	(9122.441,3946.407)
	(9133.926,3987.906)
	(9145.267,4028.479)
	(9156.467,4068.132)
	(9167.532,4106.873)
	(9189.264,4181.649)
	(9210.491,4252.864)
	(9231.240,4320.577)
	(9251.539,4384.847)
	(9271.415,4445.732)
	(9290.895,4503.290)
	(9310.007,4557.579)
	(9328.779,4608.658)
	(9347.238,4656.586)
	(9365.411,4701.420)
	(9383.326,4743.219)
	(9401.010,4782.041)
	(9435.796,4850.990)
	(9469.990,4908.732)
	(9503.810,4955.734)
	(9537.476,4992.465)
	(9605.227,5036.976)
	(9675.000,5046.000)

\path(9675,5046)	(9748.058,4984.328)
	(9774.493,4919.594)
	(9785.610,4880.330)
	(9795.488,4837.068)
	(9804.247,4790.260)
	(9812.009,4740.354)
	(9818.894,4687.804)
	(9825.023,4633.057)
	(9830.517,4576.566)
	(9835.497,4518.781)
	(9840.084,4460.153)
	(9844.399,4401.131)
	(9848.562,4342.167)
	(9852.694,4283.711)
	(9856.916,4226.214)
	(9861.350,4170.126)
	(9866.115,4115.898)
	(9871.333,4063.981)
	(9877.124,4014.824)
	(9883.610,3968.879)
	(9890.911,3926.595)
	(9899.149,3888.425)
	(9918.915,3826.224)
	(9975.000,3771.000)

\path(9975,3771)	(10025.558,3785.301)
	(10072.866,3833.580)
	(10095.645,3871.630)
	(10118.024,3919.575)
	(10140.139,3977.882)
	(10162.129,4047.019)
	(10173.119,4085.794)
	(10184.130,4127.452)
	(10195.178,4172.050)
	(10206.280,4219.648)
	(10217.454,4270.303)
	(10228.717,4324.074)
	(10240.085,4381.020)
	(10251.577,4441.198)
	(10263.209,4504.667)
	(10274.998,4571.485)
	(10286.962,4641.711)
	(10299.117,4715.403)
	(10305.272,4753.567)
	(10311.482,4792.620)
	(10317.748,4832.568)
	(10324.072,4873.419)
	(10330.458,4915.181)
	(10336.906,4957.860)
	(10343.419,5001.464)
	(10350.000,5046.000)

\path(4200,12921)	(4240.636,12881.993)
	(4279.936,12844.450)
	(4354.679,12773.635)
	(4424.531,12708.309)
	(4489.794,12648.223)
	(4550.770,12593.130)
	(4607.762,12542.785)
	(4661.072,12496.938)
	(4711.001,12455.344)
	(4757.852,12417.754)
	(4801.927,12383.922)
	(4843.529,12353.601)
	(4882.958,12326.543)
	(4956.510,12281.228)
	(5025.000,12246.000)

\path(5025,12246)	(5063.159,12228.094)
	(5107.050,12207.302)
	(5155.826,12184.486)
	(5208.643,12160.509)
	(5264.656,12136.232)
	(5323.020,12112.518)
	(5382.888,12090.229)
	(5443.417,12070.226)
	(5503.762,12053.373)
	(5563.076,12040.531)
	(5620.515,12032.563)
	(5675.233,12030.330)
	(5726.386,12034.696)
	(5773.128,12046.521)
	(5814.615,12066.668)
	(5850.000,12096.000)

\path(5850,12096)	(5879.378,12162.562)
	(5881.163,12204.179)
	(5874.810,12248.033)
	(5860.669,12291.632)
	(5839.087,12332.486)
	(5775.000,12396.000)

\path(5775,12396)	(5716.131,12424.781)
	(5655.798,12442.307)
	(5594.390,12449.767)
	(5532.294,12448.348)
	(5469.900,12439.237)
	(5407.595,12423.621)
	(5345.770,12402.690)
	(5284.811,12377.629)
	(5225.109,12349.626)
	(5167.051,12319.870)
	(5111.026,12289.547)
	(5057.423,12259.845)
	(5006.630,12231.952)
	(4959.036,12207.055)
	(4915.030,12186.342)
	(4875.000,12171.000)

\path(4875,12171)	(4800.684,12136.037)
	(4758.083,12110.208)
	(4712.541,12080.471)
	(4664.577,12048.068)
	(4614.708,12014.241)
	(4563.450,11980.233)
	(4511.321,11947.286)
	(4458.839,11916.643)
	(4406.520,11889.546)
	(4354.881,11867.237)
	(4304.441,11850.959)
	(4255.715,11841.953)
	(4209.222,11841.464)
	(4165.478,11850.732)
	(4125.000,11871.000)

\path(4125,11871)	(4081.429,11930.966)
	(4070.234,11972.235)
	(4066.260,12017.014)
	(4069.629,12062.291)
	(4080.461,12105.055)
	(4125.000,12171.000)

\path(4125,12171)	(4173.647,12200.907)
	(4226.344,12219.008)
	(4282.438,12226.576)
	(4341.281,12224.886)
	(4402.222,12215.211)
	(4464.609,12198.824)
	(4527.794,12176.998)
	(4591.125,12151.009)
	(4653.952,12122.128)
	(4715.625,12091.630)
	(4775.493,12060.789)
	(4832.906,12030.877)
	(4887.214,12003.169)
	(4937.766,11978.937)
	(4983.911,11959.457)
	(5025.000,11946.000)

\path(5025,11946)	(5081.461,11931.132)
	(5146.421,11913.795)
	(5182.843,11903.995)
	(5222.296,11893.330)
	(5265.082,11881.718)
	(5311.504,11869.076)
	(5361.863,11855.323)
	(5416.461,11840.376)
	(5475.601,11824.152)
	(5539.585,11806.569)
	(5608.715,11787.545)
	(5683.292,11766.997)
	(5722.718,11756.126)
	(5763.620,11744.843)
	(5806.035,11733.138)
	(5850.000,11721.000)

\path(2550,10221)	(2590.636,10181.993)
	(2629.936,10144.450)
	(2704.679,10073.635)
	(2774.531,10008.309)
	(2839.794,9948.223)
	(2900.770,9893.130)
	(2957.762,9842.785)
	(3011.072,9796.938)
	(3061.001,9755.344)
	(3107.852,9717.754)
	(3151.927,9683.922)
	(3193.529,9653.601)
	(3232.958,9626.543)
	(3306.510,9581.228)
	(3375.000,9546.000)

\path(3375,9546)	(3413.159,9528.094)
	(3457.050,9507.302)
	(3505.826,9484.486)
	(3558.643,9460.509)
	(3614.656,9436.232)
	(3673.020,9412.518)
	(3732.888,9390.229)
	(3793.417,9370.226)
	(3853.762,9353.373)
	(3913.076,9340.531)
	(3970.515,9332.563)
	(4025.233,9330.330)
	(4076.386,9334.696)
	(4123.128,9346.521)
	(4164.615,9366.668)
	(4200.000,9396.000)

\path(4200,9396)	(4229.378,9462.562)
	(4231.163,9504.179)
	(4224.810,9548.033)
	(4210.669,9591.632)
	(4189.087,9632.486)
	(4125.000,9696.000)

\path(4125,9696)	(4066.131,9724.781)
	(4005.798,9742.307)
	(3944.390,9749.767)
	(3882.294,9748.348)
	(3819.900,9739.237)
	(3757.595,9723.621)
	(3695.770,9702.690)
	(3634.811,9677.629)
	(3575.109,9649.626)
	(3517.051,9619.870)
	(3461.026,9589.547)
	(3407.423,9559.845)
	(3356.630,9531.952)
	(3309.036,9507.055)
	(3265.030,9486.342)
	(3225.000,9471.000)

\path(3225,9471)	(3150.684,9436.037)
	(3108.083,9410.208)
	(3062.541,9380.471)
	(3014.577,9348.068)
	(2964.708,9314.241)
	(2913.450,9280.233)
	(2861.321,9247.286)
	(2808.839,9216.643)
	(2756.520,9189.546)
	(2704.881,9167.237)
	(2654.441,9150.959)
	(2605.715,9141.953)
	(2559.222,9141.464)
	(2515.478,9150.732)
	(2475.000,9171.000)

\path(2475,9171)	(2431.429,9230.966)
	(2420.234,9272.235)
	(2416.260,9317.014)
	(2419.629,9362.291)
	(2430.461,9405.055)
	(2475.000,9471.000)

\path(2475,9471)	(2523.647,9500.907)
	(2576.344,9519.008)
	(2632.438,9526.576)
	(2691.281,9524.886)
	(2752.222,9515.211)
	(2814.609,9498.824)
	(2877.794,9476.998)
	(2941.125,9451.009)
	(3003.952,9422.128)
	(3065.625,9391.630)
	(3125.493,9360.789)
	(3182.906,9330.877)
	(3237.214,9303.169)
	(3287.766,9278.937)
	(3333.911,9259.457)
	(3375.000,9246.000)

\path(3375,9246)	(3431.461,9231.132)
	(3496.421,9213.795)
	(3532.843,9203.995)
	(3572.296,9193.330)
	(3615.082,9181.718)
	(3661.504,9169.076)
	(3711.863,9155.323)
	(3766.461,9140.376)
	(3825.601,9124.152)
	(3889.585,9106.569)
	(3958.715,9087.545)
	(4033.292,9066.997)
	(4072.718,9056.126)
	(4113.620,9044.843)
	(4156.035,9033.138)
	(4200.000,9021.000)

\put(3150,13896){\makebox(0,0)[lb]{\smash{{{\SetFigFont{8}{9.6}{rm}$N^{4S,1}(0,0)=225$}}}}}
\put(3150,11196){\makebox(0,0)[lb]{\smash{{{\SetFigFont{8}{9.6}{rm}$N^{3S+F,1}(0,e_1)=225$}}}}}
\put(5400,8421){\makebox(0,0)[lb]{\smash{{{\SetFigFont{8}{9.6}{rm}$N^{2S+2F,1}(0,2e_1)=20$}}}}}
\put(5250,5721){\makebox(0,0)[lb]{\smash{{{\SetFigFont{8}{9.6}{rm}$N^{2S+2F,1}(0,e_2)=30$}}}}}
\put(8850,5721){\makebox(0,0)[lb]{\smash{{{\SetFigFont{8}{9.6}{rm}$N^{2S+2F,1}(e_1,e_1) = 20$}}}}}
\put(9300,3021){\makebox(0,0)[lb]{\smash{{{\SetFigFont{8}{9.6}{rm}$N^{S+3F,0}(0,3e_1)=1$}}}}}
\put(5625,3021){\makebox(0,0)[lb]{\smash{{{\SetFigFont{8}{9.6}{rm}$N^{2S+2F,1}(2e_1,0)=17$}}}}}
\put(1725,3021){\makebox(0,0)[lb]{\smash{{{\SetFigFont{8}{9.6}{rm}$N^{2S+2F,1}(e_2,0)=15$}}}}}
\put(675,21){\makebox(0,0)[lb]{\smash{{{\SetFigFont{8}{9.6}{rm}$N^{S+3F,0}(0,e_1+e_2)=4$}}}}}
\put(4800,21){\makebox(0,0)[lb]{\smash{{{\SetFigFont{8}{9.6}{rm}$N^{S+3F,-1}(0,3e_1)=7$}}}}}
\put(0,5721){\makebox(0,0)[lb]{\smash{{{\SetFigFont{8}{9.6}{rm}$N^{2S+2F,0}(0,2e_1) = 96+9=105$}}}}}
\put(1500,8421){\makebox(0,0)[lb]{\smash{{{\SetFigFont{8}{9.6}{rm}$N^{3S+F,1}(e_1,0)=185$}}}}}
\put(5550,10521){\makebox(0,0)[lb]{\smash{{{\SetFigFont{8}{9.6}{rm}$\times 2$}}}}}
\put(4725,5046){\makebox(0,0)[lb]{\smash{{{\SetFigFont{8}{9.6}{rm}$\times 2$}}}}}
\put(9675,5271){\makebox(0,0)[lb]{\smash{{{\SetFigFont{8}{9.6}{rm}$\times 3$}}}}}
\put(3300,2571){\makebox(0,0)[lb]{\smash{{{\SetFigFont{8}{9.6}{rm}$\times 2$}}}}}
\put(5925,7446){\makebox(0,0)[lb]{\smash{{{\SetFigFont{8}{9.6}{rm}$\times
2$}}}}}

\end{picture}
}
\end{center}
\caption{Calculating $N^{4S,1}(0,0) = 225$.}
\label{ruledfig}
\end{figure}
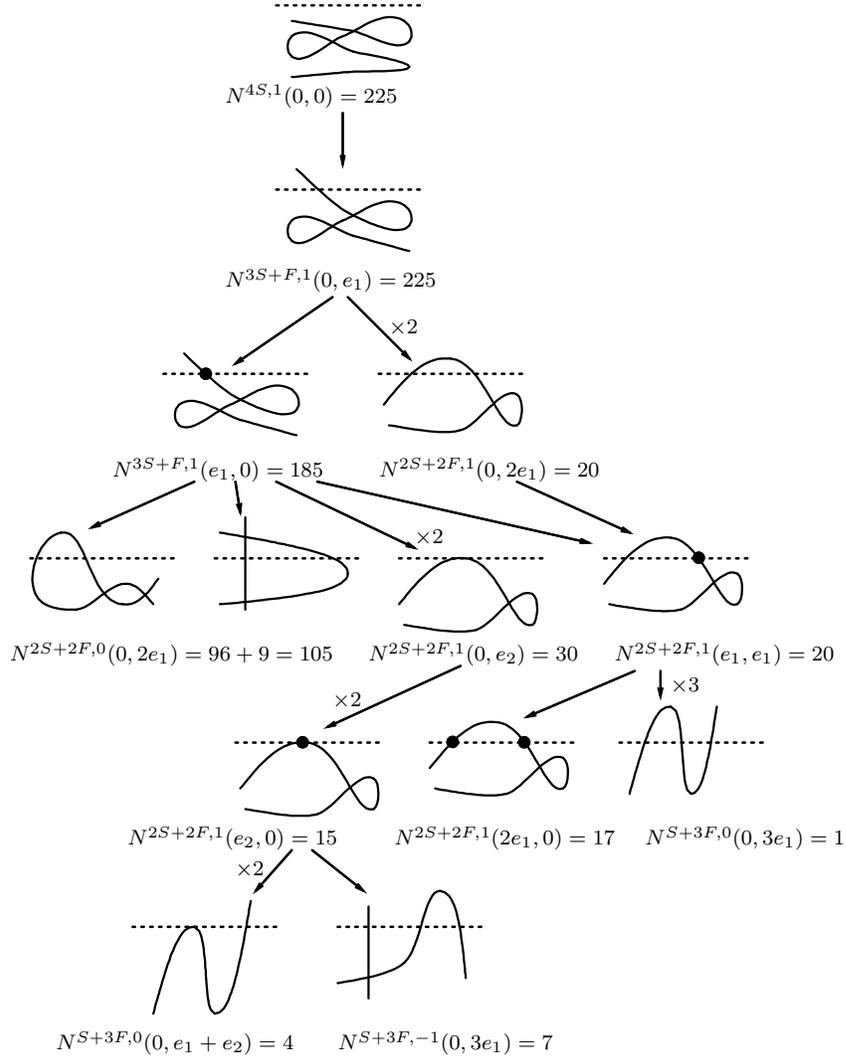

The calculation is informally summarized pictorially in Figure
\ref{ruledfig}.  The divisor $E$ is represented by the horizontal doted
line, and fixed points on $E$ are represented by fat dots.  Part of the
figure, the calculation that $N^{2S+2F,0}(0,2e_1) = 105$, has been
omitted. 

Table \ref{class2c} gives the number of genus $g$ curves in certain classes
on certain $\fn$.  Where the number of irreducible curves is different, it
is given in brackets.  Tables \ref{classf1} and \ref{classf2} give more
examples; only the total number is given, although the number of
irreducible curves could also be easily computed (using Theorem
\ref{irecursion}).  Many of these numbers were computed by a maple program
written by L. G\"{o}ttsche to implement the algorithm of Theorem
\ref{recursion}.

\begin{scriptsize}
\begin{table}
\begin{center}
\begin{tabular}{c|c|c|c|c|c|}
& $\eff_0$ & $\eff_1$ & $\eff_2$ & $\eff_3$ & $\eff_4$ \\
\hline
$2S$ & & $g=0$: 1 & $g=1$: 1 & $g=2$: 1 & $g=3$: 1 \\
     & &   & $g=0$: 10 & $g=1$: 17 & $g=2$: 24 \\
     & & & & $g=0$: 69 & $g=1$: 177 \\
& & & & & $g=0$: 406 \\

$2S+F$ & $g=0$: 1 & $g=1$: 1 & $g=2$: 1 & $g=3$: 1 \\
     &   & $g=0$: 12 & $g=1$: 20 & $g=2$: 28 \\
     & & & $g=0$: 102 (93) & $g=1$: 246 (234) \\
& & & & $g=0$: 781 (594) \\

$2S+2F$ & $g=1$: 1 & $g=2$: 1 & $g=3$: 1 \\
        & $g=0$: 12 & $g=1$: 20 & $g=2$: 28 \\
      & & $g=0$: 105 (96) & $g=1$: 252 (240) \\
& &  & $g=0$: 856 (636) \\

$2S+3F$  & $g=2$: 1 & $g=3$: 1 \\
         & $g=1$: 20 & $g=2$: 28 \\
      &  $g=0$: 105 (96) & $g=1$: 252 (240) \\
&  & $g=0$: 860 (640) \\

$2S+4F$   & $g=3$: 1 \\
          & $g=2$: 28 \\
      &  $g=1$: 252 (240) \\
&  $g=0$: 860 (640) \\
\end{tabular}
\end{center}
\caption{Number of genus $g$ curves in class $2S+kF$ on $\fn$}
\label{class2c}
\end{table}
\end{scriptsize}

\begin{table}
\begin{center}
\begin{tabular}{|c|c|c|}
\hline
Class & Genus & Number \\
\hline
\hline
$3S$ & -2 & 15 \\
   & -1 & 21 \\
   & 0 & 12 \\
   & 1 & 1 \\
\hline
$3S+F$ & 0 & 675 \\
& 1 & 225 \\
& 2 & 27 \\
& 3 & 1 \\
\hline
$3S+2F$ & 0 & 22647 \\
& 1 & 14204 \\
& 2 & 4249 \\
& 3 & 615 \\
& 4 & 41 \\
& 5 & 1 \\
\hline
$3S+3F$ & 0 & 642434 \\
& 1 & 577430 \\
& 2 & 291612 \\
& 3 & 83057 \\
& 4 & 13405 \\
& 5 & 1200 \\
& 6 & 55 \\
& 7 & 1 \\
\hline
\end{tabular}
\end{center}
\caption{Number of genus $g$ curves in various classes on $\eff_1$}
\label{classf1}
\end{table}

\begin{table}
\begin{center}
\begin{tabular}{|c|c|}
\hline
Genus & Number \\
\hline
\hline
-2 & 280 \\
-1 & 1200 \\
0 & 2397 \\
1 & 1440 \\
2 & 340 \\
3 & 32 \\
4 & 1 \\
\hline
\end{tabular}
\end{center}
\caption{Number of genus $g$ curves in class $3S$ on $\eff_2$}
\label{classf2}
\end{table}

We next review some earlier results.  (This is only a partial summary of
the voluminous research done on the subject.)  In each case but the last
two, the numbers have been checked to agree with those produced by the
algorithm given here for ``small values''.  Many of these verifications
have been done by D. Watabe, and the author is grateful to him for this.
For example, he has verified that the formula for the number of rational
curves in the divisor class $2S$ on $\fn$ passing through the appropriate
number of general points for $n \leq 9$ agrees with the numbers obtained by
Caporaso and Harris.  In each case it seems difficult to directly prove
that the numbers will always be the same, other than by noting that they
count the same thing.
\begin{itemize}
\item
The number of degree $d$ genus $g$ plane curves through $3d+g-1$ fixed
general points was calculated by Ran (cf. [R]) and Caporaso and Harris
(cf. [CH3]).  This is also the number of genus $g$ curves in the divisor
class $dS$ on $\eff_1$ through $3d+g-1$ fixed general points, and the
number of genus $g$ curves in the divisor 
class $(d-1)S+F$ on $\eff_1$ through $3d+g-2$ fixed general points. 
\item The surfaces $\eff_0$ and $\eff_1$ are convex, so the ideas of [KM]
allow one to count (irreducible) rational curves in all divisor classes on
these surfaces (see [DI] for further discussion).  These are known as the
genus 0 Gromov-Witten invariants of $\eff_0$ and $\eff_1$.
\item D. Abramovich and A. Bertram have proved several (unpublished)
formulas counting (irreducible) rational curves in certain classes on
$\eff_n$ ([AB]).  If $N^g_{\eff_n}(aS+bF)$ is the number of irreducible
genus $g$ curves in class $aS+bF$ through the appropriate number of points,
then they have shown:
\begin{enumerate}
\item[(AB1)]$N^0_{\eff_0}(aS+(b+a)F) = \sum_{i=0}^{a-1} \binom {b+2i} i
N^0_{\eff_2}(aS+bF-iE).$
\item[(AB2)]$N^0_{\eff_n}(2S+bF) = N^0_{\eff_{n-2}}(2S + (b+2)F) -
\sum_{l=1}^{n-1} \binom {2 (n+b) + 3} {n-l-1} \left( l^2 (b+2) + \binom l 2
\right).$
\item[(AB3)]  $N^0_{\eff_n}(2S) = 2^{2n} (n+3) - (2n+3) \binom {2n+1} n.$
\item[(AB4)]  $N^0_{\eff_n}(2S+bF) = N^0_{\eff_{n-1}}(2S+(b+1)F) -
\sum_{l=1}^{n-1} \binom {2(n+b) + 2} {n-l-1} l^2 (b+2).$
\end{enumerate}

Their method for (AB1) and (AB2) is to study the moduli space of stable maps of
rational curves to
$\eff_n$, and deform the surface to $\eff_{n-2}$.  For (AB3) and (AB4), they
relate curves on $\eff_n$ to curves on $\eff_{n-1}$.
\item By undoubtedly similar methods, the author has obtained the formula
\begin{eqnarray*}
N^g_{\eff_n}(2 S + k F) &=& N^g_{\eff_{n+1}}(2S + (k-1) F) \\ 
& & +
\sum_{f=0}^{n-g-1} \sum \binom {\al_1} {|\al| - g- 1} \binom { |\al|} {
\al_1, \dots, \al_n} \binom { 2n+2k+2+g} f I^{2 \al}
\end{eqnarray*}
where the second sum is over all integers $f$ and
sequences $\al$ such that $I \al = n+k-f$, $|\al| = k+1+g$, $k \leq \al_1$.
\item Caporaso and Harris (in [CH1] and [CH2]) obtained recursive formulas
for $N^0_{\eff_n}(aS+bF)$ when $n \leq 3$, and the remarkable result that
$N^0_{\eff_n}(2S)$ is the co-efficient of $t^n$ in $(1+t)^{2n+3}/ (1-t)^3$.
\item Coventry   has also recently derived a recursive formula for
the number of rational curves in {\em any} class in $\eff_n$ ([Co]), using a
generalization of the ``rational fibration'' method of [CH2].
\item 
Kleiman and Piene have examined systems with an arbitrary, but fixed,
number $\de$ of nodes ([K2]).  The postulated number of $\de$-nodal curves
is given (conjecturally) by a polynomial, and they determine the leading
coefficients, which are polynomials in $\de$.  Vainsencher determined the
entire polynomial for $\de \leq 6$ ([Va]).  Kleiman and Piene extended his
work so that a refinement of his computation of the polynomial for a given
$\de$ ought to yield the coefficients in the top $2\de$ total degrees.
They have done this explicitly for the plane for $\de \leq 4$, in
particular supporting and extending the conjecture on p. 86 of [DI], and
expect to get it for $\de \leq 6$.
\end{itemize}

\section{Dimension counts}

In this section, we prove the main dimension count we need:

\begin{tm}
\label{bigdim}
\begin{enumerate}
\item[(a)]  Each component of $\vdg_m\abG$ is of dimension at most 
$$
\Up = \Up^{D,g}(\be) = -(K_X + E) \cdot D + | \be | + g-1.
$$
\item[(b)]
The stable map $(C, \pi)$ corresponding to a general point of any
component of dimension $\Up$ satisfies the following properties.
\begin{enumerate}
\item[(i)] The curve $C$ is smooth, and the map $\pi$ is an immersion.
\item[(ii)] The image is a reduced curve.  If the $\{ p_{i,j} \}$ are
distinct, the image is smooth along its intersection with
$E$.
\end{enumerate}
\item[(c)]  Conversely, any component whose general map satisfies 
property (i) has dimension $\Up$.
\end{enumerate}
\end{tm}
By ``the image is a reduced curve'', we mean $\pi_*[C]$ is a sum of
distinct irreducible divisors on $X$.

Proposition \ref{idim} follows directly from Theorem \ref{bigdim}.

\begin{lm}[Arbarello-Cornalba, Caporaso-Harris]
\label{chlemma}
Let $V$ be an irreducible subvariety of the moduli space $\mbar_g(Y,\be)'$ where
$Y$ is smooth, such that if $(C,\pi)$ corresponds to the general point of $V$
then $C$ is smooth $\pi$ is birational.  Let $N = \coker(T_C
\rightarrow \pi^* T_Y)$, and let $\Nt$ be the torsion subsheaf of $N$.  Then:
\begin{enumerate}
\item[(a)] $\dim V \leq h^0(C, N / \Nt)$.
\item[(b)] Assume further that $Y$ is a surface.  Fix a smooth curve $G$ in
$Y$ and smooth points $\{ p_{i,j} \}$ of $G$, and assume that 
$$
\pi^* G = \sum_{i,j} i q_{i,j} + \sum_{i,j} i r_{i,j}
$$
with $\pi(q_{i,j}) = p_{i,j}$.  Then 
\begin{eqnarray*}
\dim V &\leq& h^0 (C, N / \Nt (- \sum_{i,j} i q_{i,j} - \sum_{i,j} (i-1)
r_{i,j})) \\
&=& h^0(C, N / \Nt (- \pi^* G + \sum_{i,j} r_{i,j})).
\end{eqnarray*}
\end{enumerate}
\end{lm}
This lemma appears (in a different guise) in Subsection 2.2 of [CH3]: (a)
is contained in Corollary 2.4 and part (b) is Lemma 2.6.  Part (a) was
proven earlier by E. Arbarello and M. Cornalba in [AC], Section 6.
Caporaso and Harris express (a) informally as: ``the first-order
deformation of the map $\pi$ corresponding to a torsion section of $N$ can
never be equisingular.''  Arbarello and Cornalba's version is slightly
stronger: ``the first-order deformation of the map $\pi$ corresponding to a
torsion section of $N$ can never preserve both the order and type of the
singularities of the image.''  

\begin{lm}
\label{dimbd}
Let $V$ be a component of $V_m^{D,g}\abG$ whose general point
corresponds to a map $\pi: C \rightarrow X$ where $C$ is a smooth
curve.  Then $\dim V \leq \Up$.  If $\pi$ is not an immersion then the
inequality is strict.
\end{lm}
\bpf Note that by the definition of $V^{D,g}_m\abG$, $\pi$ is a birational
map from $C$ to its image in $X$, so we may invoke Lemma 2.2.  The map $T_C
\rightarrow \pi^* T_X$ is injective (as it is generically injective, and
there are no nontrivial torsion subsheaves of invertible sheaves).  If $N$
is the normal sheaf of $\pi$, then the sequence
$$
0 \rightarrow T_C \rightarrow \pi^* T_X \rightarrow N \rightarrow 0
$$
is exact.  Let 
$\Nt$ be the torsion subsheaf of $N$.  The map $\pi$ is an
immersion if and only if $\Nt = 0$.  Now 
\begin{eqnarray*}
(\det N) (- \pi^* E + \sum r_{i,j})
&=& \oh_C( - \pi^* K_X + K_C - \pi^* E + \sum r_{i,j} ) \\
&=& \oh_C(- \pi^* K_X + K_C - \pi^* E + \sum r_{i,j} )\\
\end{eqnarray*}
By property P3, the divisor $-\pi^*(K_X+E) + \sum r_{i,j}$ is positive
on each component of $C$, so by Kodaira vanishing or Serre duality
$$
H^1(C, (\det N)(-\pi^* E+\sum r_{i,j})) = 0.
$$

As $N / \Nt$ is a subsheaf of $\det N$,
\begin{eqnarray}
\label{fred}
\lefteqn{h^0(C, N / \Nt (- \pi^* E + \sum r_{i,j}))} \\ 
&\leq& h^0 (C, (\det N) (- \pi^* E + \sum r_{i,j})) \nonumber \\
&=& \chi  (C, (\det N) (- \pi^* E + \sum r_{i,j})) \nonumber \\
&=& \deg N - E \cdot D + |\be| - g + 1 \nonumber \\
&=& - (K_X+E) \cdot D + | \be | + g - 1 \nonumber \\
&=& \Up. \nonumber
\end{eqnarray}
If $C'$ is a component of $C$ with $- \pi^*(K_X + E) \cdot C' = 1$,
then $\pi: C' \rightarrow X$ is an immersion by property P4.  Thus if
$\Nt \neq 0$, then it is non-zero when restricted to some component
$C''$ for which $-\pi^*(K_X + E) \cdot C'' \geq 2$.  Let $p$ be a
point on $C''$ in the support of $\Nt$.  Then $- \pi^*(K_X + E) + \sum r_{i,j} - p$ is positive
on each component of $C$, so by the same argument as above, $N/\Nt$ is
a subsheaf of $(\det N)(-p)$, so
\begin{eqnarray*}
\lefteqn{h^0(C, N/\Nt ( - \pi^* E + \sum r_{i,j} - p)) } \\
&\leq& h^0( C, (\det N) ( - \pi^* E + \sum r_{i,j} - p)) \\
&=& \Up - 1.
\end{eqnarray*}

Therefore, equality holds at (\ref{fred}) only if $\Nt = 0$, i.e.
$\pi$ is an immersion.  By Lemma \ref{chlemma}(a), the result follows.
\epf

{\noindent {\em Proof of Theorem \ref{bigdim}. }} Let $V$ be a
component of $\vdg_m \abG$ of dimension at least $\Up$, and
let $\pi: C \rightarrow X$ be the map corresponding to a general
point of $V$.  Let the normalizations of the components of $C$ be
$C(1)$, $C(2)$, \dots, $C(s)$, so $p_a( \coprod_k C(k)) \leq p_a(C)$
with equality if and only if $C$ is smooth.  Let $\be = \sum_{k=0}^s
\be(k)$ be the partition of $\be$ induced by $C = \cup_{k=1}^s C(k)$,
let $g(k)$ be the genus of $C(k)$, and let
$$
\Up(k) = (K_X + E) \cdot \pi_*[C(k)] + | \be(k) | + g(k) - 1.
$$
By the definition of $\vdg_m\abG$, $\pi|_{C(k)}$ is birational.

By Lemma \ref{dimbd}, $C$ moves in a family of dimension at most 
\begin{eqnarray}
\nonumber
\sum_{k=1}^s \Up(k) &=& \sum_{k=1}^s \left(- (K_X + E) \cdot \pi_*[C(k)] +
 | \be(k) | + g(k) - 1 \right) \\
\nonumber
&=& -(K_X+E) \cdot D + |\be| + p_a \left( \coprod C(k) \right) - 1 \\
&\leq& -(K_X+E) \cdot D + |\be| + p_a(C) - 1.
\label{eqna}
\end{eqnarray}
This proves part (a).  

If $\dim V = \Up$, then equality must hold in
(\ref{eqna}), so $C$ is smooth, and by Lemma \ref{dimbd}, $\pi$ is an
immersion.

We next prove that the image is smooth along its intersection with $E$.
 Requiring $\pi(r_{i_0,j_0})$ to be a fixed point $p$ imposes an
additional condition on $\vdg_m\ab$, as the locus of such maps
forms a variety of the form $\vdg_m(\al+ e_{i_0}, \be - e_{i_0}, \{
p_{i,j} \} \cup \{ p \})$ which has dimension 
$$
-(K_X+E) \cdot D + g-1 + ( |\be| - 1) < \dim \vdg_m(\al,\be, \{
p_{i,j} \})
$$ 
by part (a).

Thus if $\pi:  C \rightarrow X$ is the map corresponding to a general
point of $\vdg_m \ab$ then 
\begin{enumerate}
\item[(i)] $\pi(q_{i,j}) \neq \pi(q_{i',j'})$ for $(i,j) \neq (i',j')$ as
$p_{i,j} \neq p_{i',j'}$,
\item[(ii)] $\pi(r_{i,j}) \neq p_{i',j'}$ as requiring $\pi(r_{i,j})$ to be
fixed imposes a nontrivial condition on $\vdg_m\ab$,
\item[(iii)] $\pi(r_{i,j}) \neq \pi(r_{i',j'})$ for $(i,j) \neq (i',j')$ as
requiring both $r_{i,j}$ and $r_{i',j'}$ to be fixed imposes {\em two}
independent conditions on $\vdg_m\ab$.
\end{enumerate}
Thus the image is smooth along its intersection with $E$.

If a component of the image curve is nonreduced (with underlying
reduced divisor $D_1$), then this component cannot intersect $E$ (as
the image is smooth along $E$).  As each $C(k)$ is birational to
its image, the image curve must be the image of two components $C(k_1)$
and $C(k_2)$ for which $\Up(k_1) = \Up(k_2) = 0$.  Then
\begin{eqnarray*}
0 &=& \Up(k_1) \\
&=& \sum_{k=1}^s (- (K_X + E) \cdot \pi_*[C(k_1)] +
 | \be(k_1) | + g(k_1) - 1) \\
&=& - K_X \cdot D_1 + g(k_1) - 1
\end{eqnarray*}
As $-(K_X + E)$ is positive on all effective divisors, we must have
$g(k_1)=0$ and $- (K_X+E) \cdot D_1 = 1$.  Thus $D_1$ is rational, and
by property P4, $D_1$ is smooth.  Moreover, $E \cdot D_1 = 0$, so $- K_X
\cdot D_1 = 0$ and therefore $D_1^2 = -1$.  But then $D_1$ is an
exceptional curve not intersecting $E$, contradicting property P2.
Thus the image curve is reduced, so (b) is proved.

For part (c), let $N$ be the 
normal bundle to the map $\pi$.  As $\pi^* E$ contains no 
components of $C$,
$$
N(- K_C) = \oh_C(- \pi^*(K_X + E)) \otimes \oh_C( \pi^* E).
$$
is positive on every component of $C$ by property P3 , so 
$N$ is nonspecial.  Therefore
\begin{eqnarray*}
h^0(N) &=&  - K_X \cdot D + \deg K_C - g + 1 \\
&=& - K_X \cdot D + g-1
\end{eqnarray*}
by Riemann-Roch.  Requiring the curve to remain $i$-fold tangent to
$E$ at the point $q_{i,j}$ of $C$ (where $\pi(q_{i,j})$ is required to
be the fixed point $p_{i,j}$) imposes at most $i$ independent
conditions.  Requiring the curve to remain $i$-fold tangent to $E$
at the points $r_{i,j}$ of $C$ imposes at most $(i-1)$ independent conditions.  Thus
\begin{eqnarray*}
\dim V &\geq& - K_X \cdot D + g-1 - I \al - I \be + | \be | \\
&=& - (K_X +E) \cdot D + | \be | + g-1 \\
\end{eqnarray*}
as $I \al + I \be = D \cdot E$.
\epf

Let $V$ be an irreducible subvariety of $\mbar_g(X,D)'$, and let
$\pi: C \rightarrow X$ be the map corresponding to a general
point of a component of $V$.  Assume that $\pi^* E = \sum m q_{i,j} +
\sum m r_{i,j}$ where $\pi(q_{i,j})$ is required to be a fixed point
$p_{i,j}$ of $E$ as $C$ varies.  (In particular, no component of $C$
is mapped to $E$.)  Define $\al$ by $\al_i = \# \{ q_{i,j} \}_j$, 
$\be$ by $\be_i = \# \{ r_{i,j} \}_j$, and $\Ga= \{ p_{i,j} \}$.

\begin{pr}
\label{idimbound}
The intersection dimension of $V$ is at most 
$$
- (K_X + E) \cdot D + | \be| + g-1.
$$  
If equality
holds then $V$ is a component of $\vdg_m\abG$.
\end{pr}
The main obstacle to proving this result is that the map $\pi$ may not
map components of $C$ birationally onto their image: the map $\pi$ may
collapse components or map them multiply onto their image.

\bpf
If necessary, pass to a dominant generically finite cover of $V$ 
that will allow us to distinguish components of $C$.  
(Otherwise, monodromy on $V$ may induce a nontrivial permutation 
of the components of $C$.)

For convenience, first assume that $C$ has no contracted rational or 
elliptic components.
We may replace $C$ by its normalization; this will only make the bound
worse.  (The 
map from a component of the normalization of $C$ is also a stable map.)  We 
may further assume that $C$ is irreducible, as 
$-(K_X+E) \cdot D + | \be| + g-1$ is additive.  

Suppose $C$ maps with degree $m$ to the reduced irreducible curve $D_0 \subset X$.  Then the map $\pi: C \rightarrow D_0$ factors through the
normalization $\tilde{D}$ of $D_0$.  Let $r$ be the total ramification
index of the morphism $C \rightarrow \tilde{D}$.  By Theorem
\ref{bigdim}(a),
\begin{eqnarray*}
\idim V &\leq& \dim V \\
&\leq& - (K_X + E) \cdot D_0 + | \be| + g(\tilde{D}) - 1 \\
&=& - \frac 1 m (K_X + E) \cdot \pi_* [C]  + | \be| + \frac 1 m (
g(C) - 1 - r/2) \\
&\leq& - (K_X + E) \cdot \pi_* [C] + | \be| + g(C) - 1 
\end{eqnarray*}
where we use the Riemann-Hurwitz formula for the map $C \rightarrow
\tilde{D}$ and the fact (property P3) that $-(K_X + E) \cdot D_0 > 0$.  
Equality holds only if $m=1$, so by Theorem \ref{bigdim}, equality
holds only if $V$ is a component of $\vdg_m\abG$ for some $g$, $\al$,
$\be$, $\Ga$.

If $C$ has contracted rational or elliptic components, replace $C$
with those components of its normalization that are not contracted
elliptic or rational components (which reduces the genus of $C$) and
follow the same argument.  \epf

\begin{pr}
\label{nodal}
If $X=\fn$, the $\{ p_{i,j} \}$ are distinct, and $(C,\pi)$ is a general curve in a component of 
$\vdg_m\ab$, then $\pi(C)$ has at most nodes as singularities.
\end{pr}
Warning: To prove this, we will need more than properties P1--P4.
However, this result will not be invoked later.

\bpf
By Theorem \ref{bigdim}, $\pi$ is an immersion and $\pi(C)$ is
reduced.  Thus we need only show that $\pi(C)$ has no triple points
and that no two branches are tangent to each other.

For the former, if $s$, $t$, and $u \in C$ are distinct points of $C$
with $\pi(s) = \pi(t) = \pi(u)$, it is enough to show that there is a
section of the line bundle 
\begin{eqnarray*} 
L &:=& N ( - \sum i \cdot q_{i,j} - \sum (i-1) \cdot r_{i,j} ) \\
&=& \oh_C( - \pi^*( K_X + E) + \sum r_{i,j} + K_C)
\end{eqnarray*}
vanishing at $s$ and $t$ but not at $u$, where $N$ is 
the normal bundle to the map $\pi$.  As $\pi(C)$ 
is reduced, at most one of $s$, $t$, $u$ can lie on the component of $C$ 
mapping to the fiber $F$ through $\pi(s)$.  If none of them lie on such 
a component, then $(\pi^* F- s - t - u)$ is effective, and 
$$
L(-s-t-u - K_C) = \oh_C(\pi^*(S+F) + (\pi^* F-s-t-u) + \sum r_{i,j})
$$
is effective (consider the fiber through $\pi(s)$), so by Riemann-Roch
and Kodaira vanishing,
$$
h^0(C, L(-s-t-u)) = h^0(C,L)-3.
$$

If $u$ lies on a component of $C$ mapping to $F$, and there is a 
point $r_{i,j}$ on the same component, then both $(\pi^*F-s-t)$ and $(r_{i,j} - u)$ are both effective, and the same argument holds.  

If there is no point $r_{i,j}$ on the same component as $u$, then all
sections of $L$ vanish on $u$, and it suffices to find a section of
$L$ vanishing at $t$ but not at $s$.  But $(\pi^*F-s-t)$ is effective, so by the same argument
$$
h^0(C, L(-s-t)) = h^0(C,L)-2.
$$

To show that no two branches are tangent to each other, it is enough
to show that if $s$, $t \in C$ are distinct points with
$\pi(s)=\pi(t)$, there exists a section of $L$ vanishing at $s$ but
not at $t$, which follows from a similar argument.
\epf

The following example shows that the analogue of Proposition \ref{nodal}
does not hold for every $(X,E)$ satisfying properties P1--P4.  Let
$X=\proj^2$ and $E$ be a smooth conic.  Choose six distinct points $a$,
\dots, $f$ on $E$ such that the lines $ab$, $cd$, and $ef$ meet at a point.
Then
$$
V^{D=3L, g=-2}(\al = 6 e_1, \be=0, \Ga = \{ a, \dots, f \})
$$
consists of a finite number of maps, one of which is the map sending 
three disjoint $\proj^1$'s to the lines $ab$, $cd$, and $ef$.

\section{Identifying Potential Components}
\label{ipc}
Fix $D$, $g$, $\al$, $\be$, $\Ga$, and a general point $q$ on $E$.
Throughout this section, $\Ga$ will be assume to consist of distinct
points.  (The more general case will be dealt with in [V3].)  Let
$H_q$ be the divisor on $\vdg_m\abG$ corresponding to maps whose image
contain $q$.  In this section, we will derive a list of subvarieties (which
we will call {\em potential components}) in which each component of
$H_q$ of intersection dimension $\Up-1$ appears. (In Sections~\ref{multI}
and \ref{multII}, we will see that if the $\{ p_{i,j} \}$ are distinct,
each potential component actually appears in $H_q$.)

The potential components come in two classes that
naturally arise from requiring the curve to pass through $q$.  First, one of
the ``moving tangencies'' $r_{i,j}$ could map to $q$.  We will call
such components {\it Type I potential components}.

Second, the curve could degenerate to contain $E$ as a component.  We
will call such components {\it Type II potential components}.  For
any sequences $\al'' \leq \al$, $\ga \geq 0$, and subsets $\{
p''_{i,1}, \dots, p''_{i,\al''_i} \}$ of $\{ p_{i,1}, \dots,
p_{i,\al_i} \}$, let $g'' = g + |\ga| + 1$ and $\Ga'' = \{ p''_{i,j}
\}_{1 \leq j \leq \al''_i}$.  Define the Type II component
$K(\al'',\be,\ga,\Ga'')$ as the closure in $\mbar_g(X,D)'$ of points
representing maps $\pi: C' \cup C''
\rightarrow X$ where
\begin{enumerate}
\item[K1.] the curve $C'$ maps isomorphically to $E$,
\item[K2.] the curve $C''$ is smooth, $\pi$ maps each component of $C''$ 
birationally to its image, no 
component of $C''$ maps to $E$, and there exist $|\al''|$ points
$q_{i,j} \in C''$, $j = 1$, \dots, $\al_i''$, $|\be|$ points $r_{i,j}
\in C''$, $j = 1$, \dots, $\be_i$, $|\ga|$ points $t_{i,j} \in C''$,
$j = 1$, \dots, $\ga_i$ such that $$
\pi(q_{i,j}) = p''_{i,j} \quad \text{and} \quad
(\pi|_{C''})^*(E) = \sum i \cdot q_{i,j} + 
 \sum i \cdot r_{i,j} + 
 \sum i \cdot t_{i,j},
$$
and
\item[K3.] the intersection of the curves $C'$ and $C''$ is $\{ t_{i,j} \}_{i,j}$.
\end{enumerate}

The variety $K(\al'',\be,\ga,\Ga'')$ is empty unless $I(\al''+\be+\ga)
= (D-E) \cdot E$.  The genus of $C''$ is $g''$, and there is a degree
$\binom {\be+\ga} \be$ rational map
$$
K(\al'',\be,\ga,\Ga'') \dashrightarrow V_m^{D-E,g''}(\al'',\be+\ga,\Ga'')
$$
corresponding to ``forgetting the curve $C'$''.

\begin{tm}
\label{list}
Fix $D$, $g$, $\al$, $\be$, $\Ga$, and a point $q$ on $E$ not in $\Ga$.
Let $K$ be an irreducible component of $H_q$ with 
intersection dimension $\Up  - 1$.  Then set-theoretically, either
\begin{enumerate}
\item[I.] $K$ is a component of $\vdg_m(\al+e_k,\be-e_k,\Ga')$, 
where $\Ga'$ is the same as $\Ga$ except $p'_{k,\al_{k+1}} = q$, or
\item[II.]  $K$ is a component of $K(\al'',\be,\ga,\Ga'')$ for some $\al''$, $\ga$, $\Ga''$.
\end{enumerate}
\end{tm}
\bpf Let $(C,\pi)$ be the map corresponding to a
general point of $K$.  Consider any one-parameter subvariety
$(\cc,\Pi)$ of $\vdg_m\ab$ with central fiber $(C,\pi)$ and general
fiber not in $H_q$.  Then the total space of the curve $\cc$ in the
family is a surface, so the pullback of the divisor $E$ to this family
has pure dimension 1.  The components of $\Pi^* E$ not contained in
a fiber $\cc_t$ must intersect the general fiber and thus be the sections
$q_{i,j}$ or multisections coming from the $r_{i,j}$.  Therefore
$\pi^{-1} E$ consists of components of $C$ and points that are limits of
the $q_{i,j}$ or $r_{i,j}$.  In particular:

{\bf (*)} The number of zero-dimensional components of $\pi^* E$ not
mapped to any $p_{i,j}$ is at most $\be$, and

{\bf (**)} If there are exactly $|\be|$ such components, the
multiplicities of $\pi^* E$ at these points must be given by the sequence $\be$.

{\em Case I.}  If $C$ contains no components mapping to $E$, then 
$$
\pi^* E = \sum i \cdot a_{i,j} + \sum i  \cdot b_{i,j}
$$
where $\pi( \{ a_{i,j} \}_{i,j} ) = \{ p_{i,j} \}_{i,j} \cup \{ q \}$
and the second sum is over all $i$, $1 \leq j \leq \be'_i$ for some
sequence $\be'$.  By (*), $| \be' | \leq |\be| - 1$.  Then by
Proposition \ref{idimbound},
\begin{eqnarray*}
\idim K &\leq& - (K_X + E) \cdot D + |\be'| + g-1 \\
&\leq& - (K_X + E) \cdot D + |\be|-1 + g-1 \\
&=& \Up - 1.
\end{eqnarray*}
Equality must hold, so $|\be'| = |\be| - 1$ and $K$ is a generalized
Severi variety of maps.  The set $\pi^{-1} E$ consists of $|\al| +
|\be|$ points (which is also true of $\pi_0^{-1} E$ for a general
map $(C_0,\pi_0)$ in $\vdg_m\ab$) so the multiplicities at these
points must be the same as for the general map (i.e. $\pi^* E
|_{p_{i,j}}$ has multiplicity $i$, etc.) so $K$ must be as described
in I.

{\em Case II.} If otherwise a component of $C$ maps to $E$, say $C =
C' \cup C''$ where $C' $ is the union of irreducible components of $C$
mapping to $E$ and $C''$ is the union of the remaining components.
Define $m$ by $\pi_*[C'] = m E$, so $\pi_*[C''] = D-mE$.  Let $s = \#
(C' \cap C'' )$.

Then $p_a(C') \geq 1-m$, so 
\begin{eqnarray*}
p_a(C'') &=& g - p_a(C') + 1-s \\
&\leq& g+m-s.
\end{eqnarray*}
Assume $(\pi|_{C''})^* E = \sum i \cdot a_{i,j} + \sum i \cdot b_{i,j}$
where $\pi(a_{i,j})$ are fixed points of $E$ as $C''$ varies, and the
second sum is over all $i$ and $1 \leq j \leq \be''_i$ for some
sequence $\be''$.  By (*), $|\be''| \leq |\be| + s$.  

By restricting to an open subset if necessary, the universal map may
be written $(\cc,\Pi)$ where $\cc = \cc' \cup \cc''$, $\Pi_t(\cc'_t)
\subset E$ for all $t$, and $\Pi_t(\cc''_t)$ has no component mapping
to $E$.  Let $K'$ be the family $(\cc'',\Pi|_{\cc''})$.  We apply
Proposition \ref{idimbound} to the family $K'$:
\begin{eqnarray*}
\idim K &=& \idim K' \\
&\leq& - (K_X + E) \cdot (D-mE) + |\be''| + p_a(C'') - 1 \\
&\leq& \left( -(K_X + E) \cdot D - 2 m \right) + ( |\be| + s ) + (g+m-s) -1\\
&=& \left(   -(K_X+E) \cdot D + | \be| + g-1 \right) - 1 - (m-1) \\
&=& \Up - 1 - (m-1) \\
&\leq& \Up - 1.
\end{eqnarray*}
In the third line, we used property P1: $E$ is rational, so $(K_X+E)
\cdot E = -2$.

Equality must hold, so $m=1$ and $|\be''| = |\be| + s$.  By (**), the
multiplicity of $\pi^* E$ at the $|\be|$ points of $C''$ not in $C'
\cup \pi^* p_{i,j}$ is given by the sequence $\be$.  Let $\ga$ be
the sequence given by the multiplicities of $( \pi|_{C''})^* E$ at the
$s$ points $C' \cap C''$.  Let $\{ p''_{i,j} \}$ be the subset of $\{
p_{i,j} \}$ contained in $(\pi|_{C''})^* E$.  The only possible limits
of points $\{ q_{i,j}
\}$, $\{ r_{i,j} \}$ that could be points of $(\pi|_{C''})^{-1}
p_{i_0,j_0}$ is the section $q_{i_0,j_0}$, so $(\pi|_{C''})^*
p_{i_0,j_0}$ consists of a single point $q_{i_0,j_0}$ with
multiplicity $i_0$.  

In short, $K$ is a component of $K(\al'',\be,\ga,\Ga'')$.
\epf

There are other components of the divisor $H_q$ not counted in Theorem
\ref{list}.  For example, if $X=\proj^2$, and $E$ is a line $L$,
$D=2L$, $g=0$, $\al=2e_1$, $\be=0$, then $\vdg_m\ab$ is a
three-dimensional family (generically) parametrizing conics through 2
fixed points of $L$.  One component of $H_q$ (generically) parametrizes a line
union $L$; this is a Type II potential component.  The other
(generically) parametrizes degree 2 maps from $\proj^1$ to $L$.  This
has intersection dimension 0, so it makes no enumerative contribution.

\section{Multiplicity for Type I Intersection Components}
\label{multI}

We first solve a simpler analog of the problem.  Fix points
$\Ga=(p_{i,j})_{i,j \in {\mathbb Z}, ij \leq d}$ on $E$, and let $\al$
and $\be$ be sequences of non-negative integers.  In $\proj^d = \Sym^d
E$ representing length $d$ subschemes of $E$, we have loci
$v(\al,\be,\Ga)$ corresponding to the closure of the subvariety
parametrizing $d$-tuples
$$
\sum_i^d (\sum_{j=1}^{\alpha_i} i \cdot p_{i,j} + \sum_{k=1}^{\be_i} i \cdot r_{i,j})
$$
where $(r_{i,j})$ are any points.  Then $v(\al,\be,\Ga)$ is a smooth
variety that is the image of 
$$
\prod_{i=1}^{\be_d} \Sym^{\be_d} E
= \prod_{i=1}^{\be_d} \proj^{\be_d}
$$
embedded in $\proj^d= \Sym^d E$ by the line bundle $h_1 + 2 h_2 +
... + d h_d$ where $h_i$ is the hyperplane class of $\proj^{\be_i}$. 

Let $q$ be a point on $E$ not in $\Ga$, and let $H'_q$ be the
hyperplane in $\proj^d = \Sym^d E$ corresponding to $d$-tuples
containing $q$.  Then it is straightforward to check that, as divisors on
$v(\al,\be,\Ga)$,
\begin{equation}
\label{babybabybaby}
H'_q|_{v(\al,\be,\Ga)} = \sum_{\be_k > 0} k \cdot v(\al+e_k,\be-e_k,\Gamma_k)
\end{equation}
where $\Ga_k$ is equal to $\Ga$ in all positions except $p_{k,\al_k+1}=q$.
(For example, start by observing that equality holds set-theoretically, and
then find multiplicity by making the base change $E^{\be_k} \rightarrow
\Sym^{\be_k} E$.)

Let $K_k$ be the union of Type I potential components of the form
$\vdg_m(\al+e_k, \be-e_k, \Ga')$ as described in Theorem \ref{list}.

\begin{pr}  
\label{multItm}
The multiplicity of $H_q$ along $K_k$ is $k$.
\end{pr}
\bpf Let $U$ be the open subvariety of $\vdg_m\abG$ where $\pi^* E$
contains no components of $C$.  Then there is a rational map
$$
r:  \vdg_m\abG \dashrightarrow \Sym^{D \cdot E} E
$$
that is a morphism on $U$.  Each component of $K_k$ intersects $U$,
and $r^* H'_q = H_q$ as Cartier divisors and $r^*
v(\al+e_k,\be-e_k,\Ga_k) = K_k \cap U$ as Weil divisors, so by
(\ref{babybabybaby}) the result follows.  \epf

\section{Multiplicity for Type II Intersection Components}
\label{multII}

We translate the corresponding argument in [CH3] from the language of
the Hilbert scheme to that of maps.  

\subsection{Versal deformation spaces of tacnodes} 
We first recall facts about versal deformation spaces of tacnodes.
(This background is taken from
[CH3], Section 4, and the reader is referred there for details.)
Let $(C,p)$ be an $m^{\tth}$ order tacnode, that is, a curve
singularity analytically equivalent to the origin in the plane curve
given by the equation $y (y+x^m) = 0$.  
The Jacobian ideal $\cj$ of $y^2+yx^m$ is $(2y+x^m, yx^{m-1})$, and
the monomials $1, x, \dots, x^{m-1}$, $y, xy, \dots, x^{m-2} y$ form a
basis for the vector space $\oh / \cj$ (see [A]).  We can thus
describe the versal deformation of $(C,p)$ space explicitly.

The base is an \'{e}tale neighborhood of the origin in $\aff^{2m-1}$
with co-ordinates $a_0, \dots, a_{m-2}$, and $b_0, \dots, b_{m-1}$, and
the deformation space $\cs \rightarrow \De$ is given by
$$
y^2 + yx^m + a_0 y + a_1 xy + \dots + a_{m-2} x^{m-2} y +
b_0 + b_1 x + \dots + b_{m-1} x^{m-1} = 0.
$$
Call this polynomial $f(x,y,a_0, a_1,\dots,a_{m-2},b_0,b_1,\dots,b_{m-1})$.

There are two loci in $\De$ of interest to us.  Let $\De_m \subset \De$ be
the closure of the locus representing a curve with $m$ nodes.  This is
equivalent to requiring that the discriminant of $f$, as a function of $y$,
have $m$ double roots as a function of $x$:
$$
(x^m+ a_{m-2} x^{m-2} + \dots + a_1 x + a_0)^2 - 4 (b_{m-1}x^{m-1} + \dots + b_1
x + b_0) = 0
$$
must have $m$ double roots.  Thus $\De_m$ is given by the equations
$b_0=\dots=b_{m-1}=0$; it is smooth of dimension $m-1$.  (The locus $\De_m \subset \De$
corresponds to locally reducible curves.)

Let $\De_{m-1} \subset \De$ be the closure of the locus representing a
curve with $m-1$ nodes.  This is equivalent to the discriminant being
expressible as
$$
(x^{m-1} + \la_{m-2} x^{m-2} + \dots + \la_1 x + \la_0)^2 (x^2 +
\mu_1 x + \mu_0).
$$
From this description, we can see that $\De_{m-1}$ is irreducible of
dimension $m$, smooth away from $\De_m$, with $m$ sheets of $\De_{m-1}$
crossing transversely at a general point of $\De_m$.

Let $m_1$, $m_2$, \dots be any sequence of positive integers, and
$(C_j, p_j)$ be an $(m_j)^{\tth}$ order tacnode.  Denote the versal
deformation space of $(C_j,p_j)$ by $\De_j$,
and let $(a_{j,m_j-2}, \dots, a_{j,0}, b_{j,m_j-1}, \dots, b_{j,0})$ be
coordinates on $\De_j$ as above.  For each $j$, let $\De_{j,m_j}$ and
$\De_{j,m_j-1} \subset \De_j$ be as above the closures of loci in
$\De_j$ over which the fibers of $\pi_j$ have $m_j$ and $m_j-1$ nodes
respectively.  Finally, set
$$
\De = \De_1 \times \De_2 \times \dots,
$$
$$
\De_m = \De_{1,m_1} \times \De_{2,m_2} \times \dots,
$$
$$
\De_{m-1} = \De_{1,m_1-1} \times \De_{2,m_2-1} \times \dots .
$$
Note that $\De$, $\De_m$ and $\De_{m-1}$ have dimensions $\sum (2 m_j -
1)$, $\sum (m_j - 1)$ and $\sum m_j$ respectively.

Let $W \subset \De$ be a smooth subvariety of dimension $\sum (m_j-1)
+ 1$, containing the linear space $\De_m$.  Suppose that the tangent
plane to $W$ is not contained in the union of hyperplanes $\cup_j \{
b_{j,0}=0 \} \subset \De$.  Let $\ka := \prod m_j / \lcm(m_j)$.  Then:

\begin{lm}
\label{tacnode}
With the hypotheses above, in an \'{e}tale neighborhood of the origin
in $\De$,
$$
W \cap \De_{m-1} = \De_m \cup \Ga_1 \cup \Ga_2 \cup \dots \cup \Ga_{\ka}
$$
where $\Ga_1$, \dots, $\Ga_{\ka} \subset W$ are distinct reduced unibranch
curves having intersection multiplicity exactly $\lcm(m_j)$ with $\De_m$ at
the origin.  
\end{lm}
This lemma arose in conversations with J. Harris, and appears (with
proof) as part of [CH3] Lemma 4.3.  The key ingredient is the special
case of a single tacnode, which is [CH1] Lemma 2.14.  Results of
a similar flavor appear in [V2] Section 1.

\subsection{Calculating the multiplicity}

Suppose $K = K(\al'', \be, \ga, \Ga'')$ is a Type II component of $H_q$ (on
$\vdg_m\abG$).  Assume that the $\{ p_{i,j} \}$ are distinct.  (However,
the arguments carry through without change so long as $\{ p'_{i,j} \} = \Ga
\setminus \Ga'$ are distinct; this will be useful in [V3].)  Let $m_1$,
\dots, $m_{|\ga|}$ be a set of positive integers with $j$ appearing $\ga_j$
times ($j = 1$, 2, \dots), so $\sum m_i = I \ga$.

\begin{pr}
\label{multIItm}
The multiplicity of $H_q$ along $K$ is $m_1 \dots m_{|\ga|} = I^{\ga}$.
\end{pr}
The proof of this proposition will take up the rest of this section.

Fix general points $s_1$, \dots, $s_{\Up-1}$ on $X$, and let $H_i$
be the divisor on $\vdg_m\ab$ corresponding to requiring the image
curve to pass through $s_i$.  By Kleiman-Bertini, the intersection of
$\vdg_m\ab$ with $\prod H_i$ is a curve $V$ and the intersection of
$K$ with $\prod H_i$ is a finite set of points (non-empty as $K$ has
intersection dimension $\Up -1$).  Choose a point $(C, \pi)$ of $K
\cap H_1 \cap \dots \cap H_{\Up - 1}$.  The
multiplicity of $H_q$ along $K$ on $\vdg_m\ab$ is the multiplicity of
$H_q$ at the point $(C,\pi)$ on the curve $V$.

For such $(C,\pi)$ in $K(\al'',\be,\ga,\Ga'')$ there are unique
choices of points $\{ q_{i,j} \}$, $\{ r_{i,j}\}$ on $C$ (up to
permutations of $\{ r_{i,j} \}$ for fixed $i$): if $C = C' \cup C''$
(with $C'$ mapping isomorphically to $E$ and $\pi^{-1} E$ containing
no components of $C''$), then the condition $(\pi |_{C''})^* E = \sum
i \cdot q''_{i,j} + \sum i \cdot r_{i,j}$ with $\pi(q_{i,j}'') =
p''_{i,j}$ specifies the points $\{ q''_{i,j} \}$, $\{ r_{i,j} \}$,
and $q'_{i,j} = (\pi|_{C'})^{-1} p'_{i,j}$ specifies the points 
$\{ q'_{i,j} \} = \{ q_{i,j} \} \setminus \{ q''_{i,j} \}$.

Define the map $(\tc, \tpi)$ as follows: $C \stackrel{\pi}{\rightarrow}
X$ factors through
$$
C \stackrel{\nu}{\rightarrow} \tc \stackrel{\tpi}{\rightarrow} X.
$$
where $\nu$ is a homeomorphism and $\tpi$ is an immersion.  Each node of
$C$ is mapped to a tacnode (of some order) of $\tc$, and $\nu: C
\rightarrow \tc$ is a partial normalization.  Then $\tc$ has
arithmetic genus $\tg := g + \sum (m_i - 1)$.

Let $\Def(\tc,\tpi)$ be the deformations of $(\tc,\tpi)$ preserving
the incidences through $s_1$, \dots, $s_{\Up-1}$ and the tangencies
($\tpi^* E = \sum i \cdot q_{i,j} + \sum i \cdot r_{i,j}$,
$\tpi(q_{i,j}) = p_{i,j}$).

\begin{lm}
The space $\Def (\tc,\tpi)$ is smooth of dimension $\sum( m_j - 1) +1$.
\end{lm}
\bpf 
We will show the equivalent result: the vector space of first-order
deformations of $(\tc, \tpi)$ preserving the tangency conditions (but
not necessarily the incidence conditions $s_1$, \dots, $s_{\Up- 1}$)
has dimension $\Up + \sum( m_i - 1)$, and they are unobstructed.

As $(\tc, \tpi)$ is an immersion, there is a normal bundle to $\tpi$
$$
N_{\tc /
X} = \oh_{\tc}( - \tpi^* K_X + K_{\tc}).
$$
By property P3, as $\tpi^*(
K_X + E - \sum r_{i,j})$ is negative on every component of
$\tc$,
\begin{equation}
\label{unobstructed}
h^1( \tc, N_{\tc / X} ( - \sum i \cdot q_{i,j} - \sum (i-1) \cdot r_{i,j})) = 0
\end{equation}
so
\begin{eqnarray*}
& & h^0( \tc, N_{\tc / X} ( - \sum i \cdot q_{i,j} - \sum (i-1) \cdot r_{i,j})) \\
&=& \chi( \tc, N_{\tc / X} ( - \sum i \cdot q_{i,j} - \sum (i-1) \cdot r_{i,j})) \\
&=& \deg( \tpi^*(-K_X - E + \sum r_{i,j})) + \deg K_{\tc} - \tg + 1 \\
&=& - (K_X + E) \cdot D + | \be | + \tg - 1 \\
&=& - (K_X + E) \cdot D + | \be | + g + \sum (m_i - 1) - 1 \\
&=& \Up + \sum ( m_i - 1).
\end{eqnarray*}
Thus there are $\Up + \sum (m_i - 1)$ first-order deformations, and by
(\ref{unobstructed}) they are unobstructed. 
\epf

For convenience, let $N := N_{\tc/X} (- \sum i \cdot q_{i,j} - \sum
(i-1) \cdot r_{i,j})$.  By the proof of the above lemma, $H^0( \tc,
N)$ is naturally the tangent space to $\Def(\tc, \tpi)$.  Now $-K_X$
restricted to $C'$ has degree $K_X \cdot E = 2 + E^2$; $K_{\tc}$
restricted to $C'$ has degree $I \ga - 2$, which is $(\deg K_{C'})$ plus
the length of the scheme-theoretic intersection of $C'$ and $C''$; and
$- \sum i \cdot q'_{i,j}$ has degree $I \al'$.  Therefore 
\begin{eqnarray*}
\deg N|_{C'} &=& 2+ E^2 + I \ga - 2 - I \al' \\
&=& D \cdot E - (D-E) \cdot E + I \ga - I \al' \\
&=& (I \al + I \be) - (I \al'' + I \be + I \ga) + I \ga - I \al' \\
&=& 0
\end{eqnarray*}
so the restriction of $N$ to $C'$ is the trivial line bundle.

Also, if $q$ is a general point on $C'$ then $h^0(\tc,N(-p)) =
h^0(\tc,N) - 1$.  ({\em Proof:} From (\ref{unobstructed}), $h^1(\tc,N)
= 0$.  By the same argument, as $\deg (K_X+E)|_E = -2$, $\tpi^* ( K_X +
E - \sum r_{i,j} + q)$ is negative on every component of $\tc$, so
$h^1(\tc, N(-p)) = 0$.  Thus $h^0(\tc,N(-p))
- h^0(\tc,N) = \chi(\tc,N(-p)) - \chi(\tc,N) = -1$.)  Thus there is a
section of $N$ that is nonzero on $C'$.
 
Let $\cj$ be the Jacobian ideal of $\tc$.  In an \'{e}tale
neighborhood of the $(C,\pi)$, there are natural maps
$$
V \stackrel{\rho}{\rightarrow} 
\Def (\tc,\tpi) \stackrel{\sigma}{\rightarrow}
\De
$$
where the differential of $\si$ is given by the natural map
\begin{equation}
\label{differential}
H^0(\tc,N) \rightarrow H^0(\tc,N \otimes (\oh_{\tc} / \cj)).
\end{equation}

\begin{lm}
In a neighborhood of the origin, the morphism 
$$
\sigma: \Def(\tc, \tpi)
\rightarrow \De
$$
is an immersion, and the tangent space to $\sigma (
\Def(\tc, \tpi))$ contains $\De_m$ and is not contained in the union
of hyperplanes $\cup_j \{ b_{j,0} = 0 \}$.
\end{lm}
\bpf 
From (\ref{differential}), the Zariski tangent space to the
divisor $\si^*( b_{j,0} = 0)$ is a subspace $Z$ of $H^0(\tc,N)$
vanishing at a point of $C'$ (the $j^{\tth}$ tacnode).  But $N |_{C'}$ is
a trivial bundle, so this subspace of sections $Z$ must vanish on all
of $C'$.  As there is a section of $N$ that is non-zero on $C'$, $Z$ has
dimension at most $h^0(\tc,N) - 1 = \dim \Def(\tc,\pi) - 1$.
This proves that $\si$ is an immersion, and that the tangent space to
$\si(\Def(\tc,\tpi))$ is not contained in $\{ b_{j,0} = 0 \}$.  

Finally, if $S$ is the divisor (on $\Def(\tc,\tpi)$) corresponding to
requiring the image curve to pass through
 a fixed general point of
$E$, then $\si(S) \subset \De_m$, as the image curve must be
reducible.  As $\si$ is an immersion, 
\begin{eqnarray}
\nonumber \sum ( m_i - 1) &=& \dim \Def(\tc,\tpi) - 1 \\
\nonumber &=& \dim S \\
\nonumber &=& \dim \si(S) \\
\label{bob} &\leq& \dim \De_m \\
\nonumber &=& \sum (m_i - 1)
\end{eqnarray}
so we must have equality at (\ref{bob}), and the linear space $\De_m=
\si(S)$ is contained in $\si(\Def(\tc,\tpi))$, and thus in the tangent
space to $\si(\Def(\tc,\tpi))$.
  \epf

Thus the image $\sigma( \Def(\tc,\tpi))$ satisfies the hypotheses of
Lemma \ref{tacnode}, so the closure of the inverse image $\si^{-1}(
\De_{m-1} \setminus \De_m)$ will have $\prod m_i / \lcm(m_i)$ reduced
branches, each having intersection multiplicity $\lcm(m_i)$ with
$\si^{-1} ( \De_m)$ and hence with the hyperplane $H_q$.  Since in a
neighborhood of $(C, \pi)$ the variety $V$ is a curve birational with
$\rho(V) = \overline{ \si^{-1} ( \De_{m-1} \setminus \De_m)}$, we
conclude that the divisor $H_q$ contains $K(\al'', \be, \ga, \Ga'')$
with multiplicity $m_1 \cdots m_{| \ga| } = I^{\ga}$.

This completes the proof of Proposition \ref{multIItm}.  As an added
benefit, we see that $\vdg_m\ab$ has $I^{\ga} / \lcm(\ga)$ branches at
a general point of $K(\al'', \be, \ga, \Ga'')$.

\section{The Recursive Formulas}
\label{recursivesection}
We now collect what we know and derive a recursive formula for the degree
of a generalized Severi variety.  Fix $D$, $g$, $\al$, $\be$, $\Ga$ so that
$\Up>0$ (e.g. $(D,g,\be) \neq (kF, 1-k,\vec 0)$ when $X=\fn$).  Assume
throughout this section that $\Ga$ consists of distinct points.  Let $H_q$
be the divisor on $\vdg_m(\al,\be,\Ga)$ corresponding to requiring the
image to contain a general point $q$ of $E$.  The components of $H_q$ of
intersection dimension $\Up - 1$ were determined in Theorem \ref{list}, and
the multiplicities were determined in Propositions \ref{multItm} and
\ref{multIItm}:
\begin{pr}  
In the Chow ring of $\vdg_m\abG$, modulo Weil divisors of intersection
dimension less than $\Up - 1$,
$$
H_q = \sum_{\be_k>0} k \cdot \vdg_m(\al+e_k,\be-e_k,\Ga \cup \{ q \}) 
+ \sum I^{\ga} \cdot K( \al'', \be, \ga, \Ga'')
$$
where the second sum is over all $\al'' \leq \al$, $\Ga'' = \{
p''_{i,j} \}_{1 \leq j \leq \al''_i} \subset \Ga$, $\ga \geq 0$,
$I(\al'' + \be + \ga ) = (D-E) \cdot E$. 
\end{pr}
Intersect both sides of the equation with $H_q^{\Up - 1}$.  As those
dimension $\Up - 1$ classes of intersection dimension less than $\Up -
1$ are annihilated by $H_q^{\Up - 1}$, we still have equality:
\begin{eqnarray*}
N^{D,g}\abG &=& H_q^{\Up} \\
&=&
\sum_{\be_k>0} k \vdg_m(\al+e_k,\be-e_k,\Ga \cup \{ q \}) \cdot H_q^{\Up - 1}\\
& & + \sum I^{\ga} \cdot K( \al'', \be, \ga, \Ga'') \cdot H_q^{\Up - 1}.  
\end{eqnarray*}
As remarked in Section \ref{ipc}, each $K(\al'', \be, \ga, \Ga'')$
admits a degree $\binom {\be + \ga}
\be$ rational map to $V^{D-E,g''}_m(\al'', \be+\ga,\Ga'')$ 
(where $g'' = g- |\ga| + 1$) corresponding to ``forgetting the component mapping to $E$'', so
$$
K(\al'', \be, \ga, \Ga'') \cdot H_q^{\Up - 1} = 
\binom {\be + \ga} \ga N^{D-E,g''}(\al'', \be+\ga,\Ga'').
$$
Therefore
\begin{eqnarray*}
N^{D,g}\abG &=&
\sum_{\be_k>0} k \vdg_m(\al+e_k,\be-e_k,\Ga \cup \{ q \}) \cdot H_q^{\Up - 1}\\
& & + \sum I^{\ga} \cdot 
\binom {\be + \ga} \ga N^{D-E,g''}(\al'', \be+\ga,\Ga'').
\end{eqnarray*}
Using this formula inductively, one sees that $N^{D,g}\abG$ is independent
of $\Ga$ (so long as the $\{ p_{i,j} \}$ are distinct).

For each $\al''$, there are $\binom \al {\al''}$ choices of $\Ga''$
(as this is the number of ways of choosing 
$\{ p''_{i,1} , \dots, p''_{i,\al''_i} \}$ from
$\{ p_{i,1} , \dots, p_{i,\al_i} \}$).  Thus
\begin{eqnarray*}
N^{D,g}\ab &=& \sum_{\be_k > 0} k N^{D,g}(\al + e_k, \be-e_k) 
\\
& & + \sum I^{\ga} {\binom \al {\al''}} \binom {\be + \ga}{\be} 
N^{D-E,g''}(\al'',\be + \ga).
\end{eqnarray*}
Renaming variables $\al' := \al''$, $\be' := \be + \ga$, $g' := g''$,
this is Theorem \ref{recursion}:

\noindent
{\bf Theorem \ref{recursion}.} {\em
If $\dim \vdg(\al,\be)>0$, then
\begin{eqnarray*}
N^{D,g}\ab = \sum_{\be_k > 0} k N^{D,g}(\al + e_k, \be-e_k) 
\\
+ \sum I^{\be'-\be} {\binom \al {\al'}} \binom {\be'}{\be} 
N^{D-E,g'}(\al',\be')
\end{eqnarray*}
where the second sum is taken over all $\al'$, $\be'$, $g'$ satisfying
$\al' \leq \al$, $\be' \geq \be$, $g-g' = |\be'-\be| - 1$, $I \al' + I \be' =
(D-E) \cdot E$.}

The numbers $N^{D,g}\ab$ can be easily inductively calculated.  If
$\Up=0$, then a short calculation shows that $\be=0$, $D=kF$, and
$g=1-k$, so $N^{D,g}\ab$ is 1 is $\al=ke_1$ and 0 otherwise.  If
$\Up>0$, then $N^{D,g}\ab$ can be calculated using Theorem
\ref{recursion}.

\subsection{Theorem \ref{recursion} as a differential equation}
Define the generating function
$$
G = \sum_{D,g,\al,\be}  N^{D,g}\ab  v^D w^{g-1} \left( \frac {x^{\al}} {\al!} \right)
y^{\be} \left( \frac {z^{\Up}} {\Up!} \right)
$$
(where $w$ and $z$ are variables, $x=(x_1, x_2, \dots)$, $y = (y_1, y_2,
\dots)$, and $\{ v^D \}_{D \text{ effective}, D \neq E}$ is a
semigroup algebra)
Then Theorem \ref{recursion} is equivalent to the differential equation
\begin{equation}
\label{diffeq}
\frac {\partial G} {\partial z} = \left( \sum k y_k \frac {\partial}
{\partial x_k} + \frac {v^E} w e^{ \sum (x_k + k w \frac \partial {\partial y_k})}
\right) G.
\end{equation}
The corresponding observation for the plane is due to E. Getzler (cf.
[Ge] Subsection 6.3), and nothing essentially new is involved here.  The 
notation is slightly different from Getzler's; the introduction of a
variable $w$ corresponding to the arithmetic genus avoids the use of a
residue.

Define the generating function 
$$
G_{\ti} = \sum_{D,g,\al,\be} {N^{D,g}_{\ti}\ab}  v^D w^{g-1} \left( \frac {x^{\al}} {\al!} \right)
y^{\be} \left( \frac {z^{\Up}} {\Up!} \right).
$$
Then by a simple combinatorial argument,
$$
G = e^{G_{\ti}}.
$$
Substituting this into (\ref{diffeq}) yields a
differential equation satisfied by $G_{\ti}$:
\begin{equation}
\label{diffeq2}
\frac {\partial G_{\ti}} {\partial z} =  \sum k y_k \frac {\partial}
{\partial x_k} G_{\ti} + \frac {v^E} w e^{ \sum (x_k + {G}_{\ti} \mid_{y_k
\mapsto y_k + kw}) - G_{\ti} }
\end{equation}
where ${G}_{\ti} \mid_{y_k \mapsto y_k+kw}$ is the same as $G_{\ti}$ except
$y_k$ has been replaced by $(y_k + k w)$.  (Once again, this should be
compared with Getzler's formula in [Ge].)

However, $N_{\ti}^{D,g}\ab$ can also be calculated directly:

\noindent
{\bf Theorem \ref{irecursion}.}
{\em 
If $\dim W^{D,g}\ab>0$, then 
\begin{eqnarray*}
N_{\ti}^{D,g}\ab &=&  \sum_{\beta_k > 0} k N_{\ti}^{D,g}(\alpha + e_k,
\beta - e_k)
\\
& & + \sum \frac 1 \si \binom {\al} {\al^1, \dots, \al^l, \al- \sum \al^i} \binom {\Up^{D,g}(\be)-1} {\Up^{D^1,g^1}(\be^1), \dots, \Up^{D^l,g^l}(\be^l)} \\
& & \cdot  \prod_{i=1}^l \binom {\be^i} {\ga^i} I^{\be^i - \ga^i} N_{\ti}^{D^i,g^i}(\al^i,\be^i) 
\end{eqnarray*}
where the second sum runs over choices of $D^i, g^i, \alpha^i, \beta^i,
\gamma^i$ ($1 \le i \le l$), where $D^i$ is a divisor class, $g^i$ is a
non-negative integer, $\alpha^i$, $\be^i$, $\ga^i$ are sequences of
non-negative integers, $\sum D^i = D-E$, $\sum \ga^i = \be$, $\be^i \gneq
\ga^i$, and $\si$ is the number of symmetries of the set $\{
(D^i,g^i,\al^i,\be^i,\ga^i) \}_{1 \leq i \leq l}$.
}

(This recursion is necessarily that produced by the differential equation (\ref{diffeq2}).)

The proof is identical, except that rather than considering all maps,
we just consider maps from connected curves.  The Type I components
that can appear are analogous.  The Type II components consist of maps
from curves $C=C(0) \cup \dots \cup C(l)$ where $C(0)$ maps
isomorphically to $E$, and $C(i)$ intersects $C(j)$ if and only if one
of $\{i,j\}$ is 0.  (In the previous ``possibly reducible'' case, we only
required ``$C(i)$ intersects $C(j)$ only if one of $\{ i,j \}$ is 0.'')

The numbers $N^{D,g}_{\ti}\ab$ can be easily inductively calculated.  If
$\Up=0$, then a short calculation shows that $\be=0$, $D=F$, and
$g=0$, so $N^{D,g}_{\ti}\ab$ is 1 is $\al=e_1$ and 0 otherwise.  If
$\Up>0$, then $N^{D,g}_{\ti}\ab$ can be calculated using Theorem
\ref{irecursion}.

\section{Theorem \ref{irecursion} for $\eff_{\nm}$ computes genus $g$
Gromov-Witten invariants of $\eff_n$}
\label{gwenumerative}

The results of this section are likely all well-known, but the
author was unable to find them in the standard literature.    By
$\eff_{\nm}$, we mean $\eff_0$ if $n$ is even and $\eff_1$ if $n$ is
odd.  

(Genus $g$) Gromov-Witten invariants were defined by Kontsevich and Manin
([KM] Section 2).  We recall their definition, closely following the
discussion in [FP] Section 7 of the genus 0 case.

The varieties $\mbar_{g,n}(X,D)$ come equipped with $n$ morphisms $\rho_1$,
\dots, $\rho_n$ to $X$, where $\rho_i$ takes the point $[C, p_1, \dots,
p_n, \mu] \in \mbar_{g,n}(X,D)$ to the point $\mu(p_i)$ in $X$.  Given
arbitrary classes $\ga_1$, \dots, $\ga_n$ in $A^* X$, we can construct the
cohomology class
$$
\rho_1^*(\ga_1) \cup \cdots \cup \rho_n^*(\ga_n)
$$
on $\mbar_{g,n}(X,D)$, and we can evaluate its homogeneous component of the
top codimension on the virtual fundamental class, to produce a number, call
a {\em genus $g$ Gromov-Witten invariant}, that we denote by $I_{g,D}(\ga_1
\dots \ga_n)$:
$$
I_{g,D}(\ga_1 \cdots \ga_n) = \int_{\mbar_{0,n}(X,D)} \rho_1^*(\ga_1) \cup
\cdots \cup \rho_n^*(\ga_n) \cup F
$$
where $F$ is the virtual fundamental class.  If the classes $\ga_i$ are
homogeneous, this will be nonzero only if the sum of their
codimensions is the ``expected dimension'' of $\mbar_{g,n}(X,D)$.

By variations of the same arguments as in [FP] (p. 35):

(I)  If $D=0$, $I_{g,D}(\ga_1 \cdots \ga_n)$ is non-zero only if 
\begin{enumerate}
\item[i)] $g=0$ and
$n=3$, in which case it is $\int_X \ga_1 \cup \ga_2 \cup \ga_3$, or
\item[ii)] $g=1$, $n=1$, and $\ga_1$ is a divisor class, in which case it
is $(\ga_1 \cdot K_X) / 24$.  (The author is grateful to T. Graber for
pointing out this fact, which is apparently well-known.  This second case
is the only part of the argument that is not essentially identical to the
genus 0 presentation in [FP].)
\end{enumerate}

(II)  If $\ga_1 = 1 \in A^0 X$, $I_{g,D}(\ga_1 \cdots \ga_n)$ is nonzero
unless $D = 0$, $g=0$, $n=3$, in which case it is $\int_X \ga_2 \cup
\ga_3$.  

(III) If $\ga_1 \in A^1 X$ and $D \neq 0$, then by the divisorial axiom
([KM] 2.2.4 or [FP] p. 35), $I_{g,D}(\ga_1 \cdots \ga_n) = \left( \int_D
\ga_1 \right) \cdot I_{g,D}(\ga_2 \cdots \ga_n)$.  

In light of these three observations, in order to compute the genus $g$
Gromov-Witten invariants for a surface, we need only compute
$I_{g,D}(\ga_1 \cdots \ga_n)$ when each  $\ga_i$ is the class of a point.
For the remainder of this section, we assume this to be the case.

Now let $X$ be a Fano surface.  The ``expected dimension'' of
$\mbar_{g,n}(X,D)$ is $-K_X \cdot D + g - 1 + n$ ([KM] Section 2).  Thus
$I_{g,D}(\ga_1 \cdots \ga_n) = 0$ unless $n=-K_X \cdot D + g-1$, and the
only components of $\mbar_g(X,D)$ contributing to the integral will be
those with intersection dimension at least $-K_X \cdot D +g-1$.  (On any
other component, $\rho_1^*(\ga_1) \cap \cdots \cap \rho_n^*(\ga_n) = 0$:  there
are no curves in such a family passing through $n$ fixed general points.)

\begin{lm}
\label{gwlemma}
Let $X$ be a Fano surface, and let $D$ be an effective divisor class
on $X$.  Suppose that $M$ is an irreducible component of $\mbar_g(X,D)$
with general map $(C,\pi)$.  Then
$$
\idim M \leq -K_X \cdot D + g-1.
$$
If equality holds and $D \neq 0$, $\pi$ is an immersion.
\end{lm}
\bpf
Note that if $D$ is an effective divisor that is not 0 or the class of a
(-1)-curve, then $-K_X \cdot D > 1$.  

Assume first that $C$ is smooth and $\pi$ is birational from $C$ to its
image.  If $D$ is the class of a (-1)-curve, the result is
immediate, so assume otherwise.  Let $N= \coker(T_C \rightarrow \pi^* T_X)$
be the normal sheaf of $\pi$.  Let $\Nt$ be the torsion subsheaf of $N$, so
$\pi$ is an immersion if and only if $\Nt = 0$.  Let $\det N$ be the
determinant line bundle.  By [AC] Section 6 or [CH3] Lemma 2.2, $\dim M
\leq h^0(C, N/\Nt)$.  (Caporaso and Harris express this informally as:
``the first-order deformations of the birational map $\pi$ corresponding to
a torsion section of $N$ can never be equisingular.''  Arbarello and
Cornalba's version, proved earlier, is slightly stronger: ``the first-order
deformation of the birational map $\pi$ corresponding to a torsion section
of $N$ can never preserve both the order and type of the singularities of
the image.'')  As $N/\Nt$ is a subsheaf of $\det N$,
\begin{eqnarray*}
\idim M &\leq & \dim M \\
& \leq & h^0(C, N / \Nt) \\
& \leq & h^0(C, \det N) \\
& = & h^0(C, \omega_C(- \pi^* K_X)) \\
& = & \chi(C, \omega_C(- \pi^* K_X)) \\
&=& \deg(K_C - \pi^* K_X) - g + 1 \\
&=& - K_X \cdot C + g-1
\end{eqnarray*}
where equality in the fifth line comes from Kodaira vanishing or Serre
duality, as $X$ is Fano.

If $\Nt \neq 0$ and $p$ is in the support of $\Nt$, then $N/\Nt$ is
actually a subsheaf of $\det N(-p)$.  But $-K_X \cdot D > 1$ (as $D \neq
0$, and $D$ is not the class of a (-1)-curve), so $-
\pi^* K_X -p$ is positive on $C$, and the same argument gives
\begin{eqnarray*}
\idim M &\leq & \chi(C, \omega_C(- \pi^* K_X - p)) \\
&=& - K_X \cdot C + g-2.
\end{eqnarray*}
Thus the lemma is true if $C$ is smooth and $\pi$ is birational.

Assume next that $C$ is smooth, but $\pi$ is not birational.
If $D=0$, the result is immediate, so assume otherwise.
The morphism $\pi$ factors through 
$$
C \stackrel {\pi_1} \rightarrow C' \stackrel {\pi_2} \rightarrow X
$$ 
where $C'$ is a smooth curve birational to
its image under $\pi_2$.  Let $d>1$ be the degree of $\pi_1$.  If $M'$ is
any irreducible component of $\mbar_{g(C')}(X, \pi_{2*}[C'] = D/d)$ then, by
the case already proved,
\begin{eqnarray*}
\idim M' &\leq& - K_X \cdot \pi_*[C'] + g(C') - 1 \\
&=& - \frac 1 d K_X \cdot D + g(C')-1.
\end{eqnarray*}
By Riemann-Hurwitz, $g-1 \geq d(g(C') - 1)$ with equality only if $\pi_1$ is
unramified, so
\begin{eqnarray*}
\idim M & \leq & - \frac 1 d K_X \cdot D + \frac 1 d (g - 1) \\
& \leq & - K_X \cdot D + g - 1
\end{eqnarray*}
with equality iff $-K_X \cdot D + g - 1 = 0$ and $\pi_1$ is unramified.  As
$X$ is Fano and $D$ is effective and non-zero, $-K_X \cdot D + g-1 = 0$ iff
$g=0$ and $K_X \cdot D = -1$, i.e. $D$ is an exceptional curve and $C'
\cong \proj^1$.  But there are no unramified maps to $C'$, so 
$$
\idim M < - K_X \cdot D + g-1.
$$
Thus the lemma is true if $C$ is smooth.

If $C$ has irreducible components with normalizations $C_1$, \dots,
$C_{l'}$ (where $C_k$ has geometric genus $g_k$, and $\pi_*[C_k] = D_k$)
and $C_1$, \dots, $C_l$ are those components that are not contracted
rational or elliptic components, then by the smooth case above,
\begin{eqnarray*}
\idim M &\leq & \sum_{k=1}^l \left( - K_X \cdot D_k + g_k - 1 \right) \\
&=& - K_X \cdot D + \sum_{k=1}^l (g_k -1).
\end{eqnarray*}
But $\sum_{k=1}^l (g_k - 1) \leq g-1$ with equality if and only if
$l=1$, $C_1$ is smooth, and there are no contracted rational or
elliptic components.
\epf

Thus the genus $g$ Gromov-Witten invariants of $\fn$ can be computed as
follows.  By our earlier comments, we need only compute $I_{g,D}(\ga_1
\cdots \ga_n)$ where $D$ is effective and nonzero, and the $\ga_i$ are
(general) points.  Genus $g$ Gromov-Witten invariants are
deformation-invariant ([LT] Theorem 6.1), and $\fn$ degenerates to
$\eff_{n+2}$ with the classes $(E,F)$ on $\fn$ transforming to $(E+F,F)$ on
$\eff_{n+2}$ (well-known; S. Katz has suggested the reference [N],
p. 9-10).  Hence if $D$ is the class $aE + bF$ on $\fn$, $I_{g,D}(\ga_1
\cdots \ga_n)$ on $\fn$ is $I_{g,D'}(\ga_1 \cdots \ga_n)$ on
$\eff_{\nm}$, where
$$
D' = \left( a - [n/2] b \right) E + b F.
$$
By Lemma \ref{gwlemma}, this is the number of immersed genus $g$ curves in
class $D'$ through the appropriate number of points of $\eff_{\nm}$.
If $D' = E$, the number is 1.  Otherwise, the number is recursively
calculated by Theorem \ref{irecursion}.

\noindent
\address{Department of Mathematics \\ Princeton University \\
Fine Hall, Washington Road \\ Princeton, NJ 08544-1000 \\ vakil@math.princeton.edu}
\end{document}